\documentclass[]{interact}
\usepackage{amssymb}
\usepackage{amsfonts}
\usepackage{mathrsfs}
\usepackage[all]{xy}
\usepackage{multicol}
\usepackage{siunitx} 
\usepackage{enumerate} 
\usepackage{pgfplots}
\usepackage{pgfplotstable}
\usepackage{tikz,pgfplots}
\usepackage{graphicx,wrapfig}
\usepackage{wasysym} 
	\usepackage{geometry}
	\usepackage{arydshln}
\usepackage{epstopdf}
\usepackage[caption=false]{subfig}

\usepackage[numbers,sort&compress]{natbib}
\bibpunct[, ]{[}{]}{,}{n}{,}{,}
\makeatletter
\def\NAT@def@citea{\def\@citea{\NAT@separator}}
\makeatother
\theoremstyle{plain}
\newtheorem{theorem}{Theorem}[section]
\newtheorem{lemma}[theorem]{Lemma}

\newtheorem{proposition}[theorem]{Proposition}

\theoremstyle{definition}
\newtheorem{definition}[theorem]{Definition}

\theoremstyle{remark}
\newtheorem{remark}{Remark}

 \newcommand{\bigzero}{\mbox{\normalfont\Large\bfseries 0}}
\begin{document}


\title{Polar Varieties in Cayley-Klein Spaces}
\author{
\name{Fahimeh Heidari\textsuperscript{a}\thanks{ \textsuperscript{a} Email: fheidari@aut.ac.ir} and Bijan Honari\textsuperscript{b} \thanks{
\textsuperscript{b} Email: honari@aut.ac.ir}}
\affil{\textsuperscript{a,b}Department of Mathematics and Computer Science, Amirkabir University of Technology (Tehran Polytechnic),
Iran} 
}
\maketitle
\begin{abstract}
In this paper, we introduce the notion of a total polar for an arbitrary subspace of a Cayley-Klein space in an analytical framework. We show that the set of all total polars of a subspace is a Schubert variety. The notion of total polar gives a definition for a subspace to be tangent to the absolute figure of the space. By specifying tangent lines, tangent cones and then spheres are defined. This definition of the sphere does not depend on the metric of the space. It is proved that every reflection of a Cayley-Klein space, defined by two subspaces which are total polar to each other, is a motion of the space. On the other hand, each motion in a Cayley-Klein space of dimension
$n$
is a product of at most
$n+1$
reflections in point-hyperplane pairs.
\end{abstract}
\begin{amscode}
51F10, 51F20, 14M15, 51F15, 51N15
\end{amscode}
\section{Introduction}
Projective geometry is all geometry. This is the famous assertion of the great English geometer Arthur Cayley appeared in his sixth Memoir upon Quantics
\cite{cayley1859iv}
in 1859. Whereas up to this time, the projective geometry had been regarded as a part of metric geometry, Cayley considered the Euclidean plane within the projective plane and obtained the Euclidean metric only by using projective notions. In Cayley's words: ‘ A chief object of the present memoir is the establishment, upon purely descriptive principles, of the notion of the distance’. It seems a paradoxical statement, because we have no distances and no angles in projective geometry. In fact, Cayley was able to construct the distance between two points and the angle between two lines by considering a quadric, namely two imaginary points on the line at infinity, in the projective plane. Ten years later, Klein in 
\cite{klein1873sogenannte}
took up this idea of Cayley and derived distances between points and angles between lines of the hyperbolic and elliptic planes by fixing quadrics in projective plane. These two non-Euclidean spaces along with the Euclidean space and some other ones such as the Minkowskian and Galilean spaces are nowadays called Cayley-Klein spaces. A Cayley-Klein space of dimension 
$n$
is an
$n$-dimensional
real projective space with a sequence of quadrics, called the absolute figure of the space, in which each quadric is defined in the vertex of the previous one and the last quadric is non-degenerate. As the literature of this field, we refer the reader to the report 
\cite{yaglom1964projective}
by I. M. Yaglom et al, as well as the textbook 
\cite{giering2013vorlesungen}
of O. Giering. 
\par
For a non-degenerate quadric
$Q$
in an
$n$-dimensional
real projective space, that yields a Cayley-Klein space with the absolute figure
$Q,$
there exists a polarity. This polarity associates to each 
$k$-dimensional
subspace
$K$
the polar of it with respect to
$Q$
which is a subspace of dimension
$n-k-1.$
Giering has used the term total polar to refer to the polar of a subspace with respect to a quadric. In general Cayley-Klein spaces, he has defined this notion only for some special subspaces which are called regular by him. For each regular subspace of dimension 
$k$
there corresponds a unique total polar which is a subspace of dimension 
$n-k-1.$
This incomplete generalization of polarity in an arbitrary Cayley-Klein space gives rise to an incomplete definition of the notions of orthogonality or reflection in these spaces. The notion of orthogonality is defined in 
\cite{giering2013vorlesungen}
only for regular subspaces. Every reflection of a Cayley-Klein space is defined in
\cite{yaglom1964projective}
by means of two subspaces that one of them is a polar of the other. Although both of these subspaces may be considered non-regular, a clear definition of polars for non-regular subspaces does not exist there. Until now few researchers have addressed the problem of associating a polar (or total polar) to an arbitrary 
subspace in a Cayley-Klein space. Richter-Gebert in his book
\cite{richter2011perspectives}
has focused on Cayley-Klein planes. His approach is completely different from Giering. He has introduced a primal-dual pair of conics for defining the absolute figure of a Cayley-Klein plane. For each point and each line of the plane there corresponds a line and a point as a polar respectively which is not necessarily unique. He has proved that each reflection of a Cayley-Klein plane is a motion of the plane. The notion of orthogonality is also defined for any two given lines by using the dual conic. H. Struve and R. Struve in
\cite{struve2004projective, struve2008lattice}
have provided a purely synthetic framework for defining Cayley-Klein spaces. They have presented a criterion in this framework for determining when two arbitrary subspaces of a Cayley-Klein space of dimension
$n$
are polar to each other. It is proved that all polars of a $k$-dimensional
subspace have dimension
$ n-k-1.$
They have shown that every reflection of a Cayley-Klein space leaves invariant the absolute figure of the space. Also the notion of orthogonality has been defined for any two given subspaces by them. Although both definitions for polars agree with Giering's definition for regular subspaces, there are some disagreements between them which are expressed in the next section. 
\par
This paper aims to give a new definition for polars in a Cayley-Klein space in an analytical framework. The definition is in agreement with 
\cite{giering2013vorlesungen}
for regular subspaces. In section two, we give an overview of the existing definitions for polars and in the third section, we introduce our definition for this notion. The definition associates to each
$ k $-dimensional
subspace of a Cayley-Klein space of dimension 
$ n,$
a family of 
$ (n-k-1)$-dimensional 
subspaces as total polars. It is proved that this family of subspaces is a Schubert variety ensuring the existence of at least one total polar for each subspace. Also we give the necessary and sufficient conditions for the uniqueness of a total polar and get some important properties of total polars. It is shown that if
$ K$
is a subspace of
$ K^{\prime}, $
then every total polar of
$ K $
contains a total polar of
$ K^{\prime} $
and every total polar of
$ K^{\prime} $
can be extended to a total polar of 
$ K. $ 
In section four, we give the definition for a subspace to be tangent to the absolute figure of a Cayley-Klein space. Achieving a complete description of tangent lines, enables us to define tangent cones and present their equations. Using these tangent cones, we will obtain the equations of all spheres of the space without applying its metric. In section five, we will give the definition of a reflection in a Cayley-Klein space. We show that every reflection is a motion of the space and each motion of a Cayley-Klein space of dimension
$n$
is a product of at most
$n+1$
reflections in point-hyperplane pairs.
\section{A Survey of Definitions of Polars}
In this section, we give a brief overview of the existing definitions of total polars (or polars) in Cayley-Klein spaces. Consider a projective space $ P(\mathrm{V}) $
for a vector space
$ \mathrm{V}$ 
over a field
$\mathbb{F}.$
That is, the quotient of
$ \mathrm{V} \setminus \lbrace 0 \rbrace   $
by the equivalence relation that identifies a nonzero vector
$ \textit{\textbf{x}}$
with any other vector
$ \alpha \textit{\textbf{x}} $
for 
$ \alpha \neq 0.$
We denote the point
$X=[\textit{\textbf{x}}]$
of
$P(\mathrm{V})$
by 
$X(\textit{\textbf{x}})$ 
if needed. We write 
$ P^{n}(\mathbb{R}) $
and
$ P^{n}(\mathbb{C}) $
rather than
$ P(\mathbb{R}^{n+1}) $
and
$ P(\mathbb{C}^{n+1}) $
respectively. For a point
$X((x_{0}, \ldots, x_{n}))$ 
of
$ P^{n}(\mathbb{R}) $
the column matrix
$(\begin{array}{c c c }
x_{0} & \ldots & x_{n} 
\end{array})^T,$
called the coordinate matrix of
$X,$
is denoted by
$\mathsf{X}.$
Also the coordinate matrix of a hyperplane of
$ P^{n}(\mathbb{R}) $
with the equation
$a_{0}x_{0}+\ldots + a_{n}x_{n}=0$
is the row marix
$(\begin{array}{c c c }
a_{0} & \ldots & a_{n} 
\end{array}).$
A quadric
$Q$
in
$ P(\mathrm{V}) $
is equal to the set
$ \lbrace X(\textit{\textbf{x}}) \vert f(\textit{\textbf{x}},\textit{\textbf{x}})=0 \rbrace $
for a symmetric bilinear form 
$ f : \mathrm{V} \times \mathrm{V} \rightarrow \mathbb{F}$
which is not identical to the zero form. We say that the point
$ X(\textit{\textbf{x}}) $
is conjugate to a point
$ Y(\textit{\textbf{y}}) $
with respect to 
 $ Q $
if
$f(\textit{\textbf{x}},\textit{\textbf{y}})=0 $.
The polar of a subspace  
$ K $
with respect to 
$ Q, $
denoted by  
$ K^{p}, $ 
is defined to be the set of all points of 
$ P(\mathrm{V}) $
that are conjugate to every point of 
$ K .$
For the simplicity of notation, for each point 
$ X $
of 
$ P(\mathrm{V}) $
we write
$ X^{p} $ 
instead of
$ \lbrace X \rbrace^{p}. $
The vertex of
$ Q $
is defined to be the polar of 
$ P(\mathrm{V}) .$ 
The quadric
$ Q $
is called non-degenerate if its vertex is empty. If 
$ X $
is a point which is not in the vertex,   
$ X^{p} $ 
is a hyperplane of
$ P(\mathrm{V}); $
otherwise 
$ X^{p}= P(\mathrm{V}).$
In the first case, 
$ X^{p} $ 
is called the polar hyperplane of 
$ X. $
If
$ \mathrm{V}^{\mathbb{C}} $
is the complex extension of a real vector space
$ \mathrm{V}, $
$ P(\mathrm{V}) $
is embeddable into its complex extension 
$P(\mathrm{V}^{\mathbb{C}})$ 
and every quadric in 
$ P(\mathrm{V}) $
can be extended to a quadric in 
$P(\mathrm{V}^{\mathbb{C}}).$
Now we give the formal definition of a Cayley-Klein space according to
\cite{giering2013vorlesungen}.
\begin{definition}
Set
$A_{0}=P^{n}(\mathbb{R})$
and let
$ Q_{i} $
be a quadric in 
$ A_{i} $
with the vertex
$ A_{i+1} $
for 
$ 0 \leq i \leq r.$
The real projective space
$ P^{n}(\mathbb{R}) $
with the following sequence 
$$   \hat{Q_{0}} \supseteq \hat{A_{1}} \supseteq \hat{Q_{1}} \supseteq \hat{A_{2}} \supseteq \cdots \supseteq \hat{A_{r}} \supseteq \hat{Q_{r}} \supseteq \hat{A}_{r+1}= \varnothing $$
is a Cayley-Klein space of dimension
$ n $ 
where for
$ 0 \leq i \leq r,$
$ \hat{Q_{i}} $
and 
$ \hat{A}_{i+1}$
are the complex extensions of
${Q_{i}} $
and
$ {A_{i+1}}$
respectively. The above sequence is called the absolute figure of the space.
\end{definition}
Giering has defined total polars for some subspaces which are called regular by him.
\begin{definition}
Consider an
$n$-dimensional
Cayley-Klein space with the absolute figure
$$   \hat{Q_{0}} \supseteq \hat{A_{1}} \supseteq \hat{Q_{1}} \supseteq \hat{A_{2}} \supseteq \cdots \supseteq \hat{A_{r}} \supseteq \hat{Q_{r}} \supseteq \hat{A}_{r+1}= \varnothing. $$
A subspace
$ K $
of
$ P^{n}(\mathbb{R}) $
satisfying
$ K \cap A_{l}\neq \varnothing$
and
$K \cap A_{l+1}=\varnothing $ 
for an integer
$ l,$
$ 0 \leq l \leq r, $
is called a regular subspace, if 
$ K+A_{l}= P^{n}(\mathbb{R}). $
\end{definition}
In this case the total polar of 
$K$
is defined to be the polar of
$ K \cap A_{l} $
with respect to
$ Q_{l}. $
By regularity, the total polar of a
$k$-dimensional
subspace is of dimension
$n-k-1.$
Throughout this paper, we denote the polar of
$ K \cap A_{j} $
with respect to
$ Q_{j} $
by
$ (K \cap A_{j})^{p_{j}}  $
for 
$ 0 \leq j \leq r. $
In particular 
$ K^{p}=K^{p_{0}}.$
\\
For instance, the absolute figure of Euclidean plane is as follows. The first conic is a line 
$ \mathrm{L}_{\infty} $
equal to its own vertex and the second one is a non-degenerate conic
$\mathrm{Q}$
consisting of two imaginary points. The line
$ \mathrm{L}_{\infty} $
which is considered as the line at infinity in Euclidean plane, is the only non-regular line of the plane. The total polar of any other line
$L$
is the polar of the point 
$L \cap \mathrm{L}_{\infty} $
with respect to 
$\mathrm{Q}.$
Each point of the plane not at infinity is a regular point
whose total polar is 
$ \mathrm{L}_{\infty}. $
\par
Richter-Gebert in
\cite{richter2011perspectives}
has defined a polar for each point and each line of a Cayley-Klein plane. He describes the absolute figure of the plane by means of a primal-dual pair of conics. The dual conic represents a set of lines consisting of all tangent lines to the primal conic. This pair of conics is given by a pair 
$ (\mathsf{A}, \mathsf{B}) $
of real symmetric nonzero
$ 3 \times 3 $
matrices that satisfy
$ \mathsf{A} \mathsf{B}=\lambda \mathsf{I} $
for a real number
$ \lambda.  $
Let 
$(\begin{array}{c c c }
x_{0} & x_{1} & x_{n} 
\end{array})^T$
and 
$(\begin{array}{c c c }
l_{0} & l_{1} & l_{2} 
\end{array})$
be the coordinate matrices of a point
$ X $
and line
$ L$
respectively.
$ X $
and 
$ L $
are said to be polar to each other if there exist constants
$ \lambda, \ \mu \in \mathbb{R}$
such that 
$ \mathsf{A} \mathsf{X}= \lambda \mathsf{L}^{T} $
and
$\mathsf{L} \mathsf{B}= \mu \mathsf{X}^{T} .$ 
For example, matrices
$ \mathsf{A} $
and
$\mathsf{ B} $
used for constructing the absolute figure of Euclidean plane are
$$  \mathsf{A} = \begin{pmatrix}
1 & 0 & 0 \\ 
0 & 0 & 0 \\ 
0 & 0 & 0
\end{pmatrix}  \text{and} \ 
\mathsf{B}= \begin{pmatrix}
0 & 0 & 0 \\ 
0 & 1 & 0 \\ 
0 & 0 & 1
\end{pmatrix}.  $$
The matrix 
$ \mathsf{A} $
gives the point conic
$ x_{0}^{2}=0$
which is denoting the line
$ \mathrm{L}_{\infty} $
and 
$\mathsf{B} $
describes the line conic 
$ l_{1}^{2} + l_{2}^{2}=0. $
This conic consists of all lines passing through the point with coordinate matrices
$(\begin{array}{c c c }
0 & 1 & i 
\end{array})^T$
or
$(\begin{array}{c c c }
0 & 1 & -i 
\end{array})^T$
denoting the two imaginary points at infinity. In this way to each regular point and line there corresponds a unique polar which is equal to the one given by Giering. All points of the plane are polars of the line
$ \mathrm{L}_{\infty}. $
For each point
$X$
on this line there has corresponded the pencil of lines through the polar of
$X$
with respect to
$\mathrm{Q}.$
\par
H. Struve and R. Struve in
\cite{struve2004projective, struve2008lattice}
establish a correspondence between subspaces of an arbitrary
Cayley-Klein space associating at least one polar to each subspace. They have defined Cayley-Klein spaces in a synthetic framework, in which the ambient structure is a projective lattice. Let 
$( L, \wedge, \vee) $
be a projective lattice with the least and the greatest elements
$ 0$
and
$1 $
respectively. Polarities of a projective lattice are anti-automorphisms of order 2. The synthetic definition of a Cayley-Klein space introduced by them is as follows.
$ ( L, \ (([\epsilon_{0}, \epsilon_{1}],\pi_{0}), \ldots, ([\epsilon_{r}, \epsilon_{r+1}],\pi_{r}))) $
is a Cayley-Klein lattice of dimension 
$ n \geq 0,$
if the following conditions are satisfied:
\begin{enumerate}
\item
$L$
is a projective lattice of finite length
$ n+1. $
\item
$ 1= \epsilon_{0} > \epsilon_{1} > \cdots > \epsilon_{r+1}=0 $
is a chain of
$L.$
\item
$ \pi_{k}, \ 0 \leq k \leq r, $
is a polarity on the interval
$ [\epsilon_{k},\epsilon_{k+1}]=\lbrace \alpha \in L \vert \epsilon_{k} \geq \alpha \geq \epsilon_{k+1} \rbrace.  $
\end{enumerate}
They bring the notion of a polar for an arbitrary element of
$L,$
by considering for each
$ k,\ 0 \leq k \leq r, $
the projection 
$ \varphi_{k}$
which is a mapping from 
$ L $
onto
$[\epsilon_{k},\epsilon_{k+1}]$
with
$ \varphi_{k}(\alpha)=(\alpha \wedge \epsilon_{k}) \vee \epsilon_{k+1}. $
For a given
$ \alpha \in L,$ 
$ \beta \in L $ 
is called a polar of 
$ \alpha $
if 
$ \pi_{k}(\varphi_{k}(\alpha))= \varphi_{k}(\beta)$
for 
$ 0 \leq k \leq r. $
Since the set of all subspaces of a projective space is a projective lattice, it has shown that every 
$n$-dimensional
Cayley-Klein space associates to a Cayley-Klein lattice of dimension
$ n $
\cite{struve2004projective}.
The lattice associated to the Cayley-Klein space with the absolute figure 
$$   \hat{Q_{0}} \supseteq \hat{A_{1}} \supseteq \hat{Q_{1}} \supseteq \hat{A_{2}} \supseteq \cdots \supseteq \hat{A_{r}} \supseteq \hat{Q_{r}} \supseteq \hat{A}_{r+1}= \varnothing. $$
is
$ (L( P^{n}(\mathbb{R})), \ (([A_{0}, A_{1}],\pi_{0}), \ldots, ([A_{r}, A_{r+1}],\pi_{r}))), $
where
$L( P^{n}(\mathbb{R}))$
is the set of all subspaces of
$ P^{n}(\mathbb{R}) $
and
$\pi_{i}$
is a polarity of
$[A_{i}, A_{i+1}]$
mapping
$K,$
$ A_{i+1} \subseteq K \subseteq A_{i},$
to
$K^{p_{i}}$
for 
$0 \leq i \leq r.$
In this way, all polars of a subspace in a Cayley-Klein space are determined. It has proved that all polars  of a 
$k$-dimensional
subspace have dimension
$ n-k-1. $
Consider the absolute figure of Euclidean plane again in which the first conic is the line
$ \mathrm{L}_{\infty} $
and the second one is the conic
$ \mathrm{Q} $
consisting of two imaginary points. Assume that 
$ \pi_{1} $
is the polarity of
$ [\mathrm{L}_{\infty}, \varnothing] $
induced by
$ \mathrm{Q} $
and 
$ \pi_{0} $
is that of 
$ [P^{2}(\mathbb{R}), \mathrm{L}_{\infty}] $
interchanging the only two elements 
$ \mathrm{L}_{\infty}  $
and 
$P^{2}(\mathbb{R})  $
of it. With these assumptions 
$ (L( P^{2}(\mathbb{R})), ([P^{2}(\mathbb{R}), \mathrm{L}_{\infty}], \pi_{0}),( [\mathrm{L}_{\infty}, \varnothing], \pi_{1})  ) $
is the Cayley-Klein lattice associated to Euclidean plane. It is easy to see that the polar of each regular point or each regular line is identical with the one given by Giering. To a point
$X$
at infinity there have corresponded all lines through the polar of 
$X$
with respect to
$\mathrm{Q}$ 
with the exception of
$\mathrm{L}_{\infty},$
and every point of the plane is a polar of
$ \mathrm{L}_{\infty} $
except the points on it. It is in disagreement with the previous result obtained from
\cite{richter2011perspectives}
in which, all lines of the pencil are polars of 
$X$
and each point of the plane is a polar of the line
$\mathrm{L}_{\infty}.$
\section{Polar Variety}
In this section, we give our definition for total polars in a Cayley-Klein space. We will prove that the set of all total polars of a subspace is a Schubert variety. Finally, some properties of total polars are provided.
\begin{definition}
Let 
$ X $
be a point in a Cayley-Klein space of dimension
$ n $
with the absolute figure
$$   \hat{Q_{0}} \supseteq \hat{A_{1}} \supseteq \hat{Q_{1}} \supseteq \hat{A_{2}} \supseteq \cdots \supseteq \hat{A_{r}} \supseteq \hat{Q_{r}} \supseteq \hat{A}_{r+1}= \varnothing. $$
Assume that 
$ X \in A_{l}\setminus A_{l+1}$
for some 
$l,$
$ 0\leq l \leq r .$ 
A total polar of 
$ X $
with respect to the absolute figure is a hyperplane of 
$ P^{n}(\mathbb{R}) $
that contains 
$ X^{p_{l}}. $
\end{definition}
\begin{definition}
Let 
$ K $
be a 
$k$-dimensional
subspace of a Cayley-Klein space of dimension 
$ n $
with the absolute figure
$$   \hat{Q_{0}} \supseteq \hat{A_{1}} \supseteq \hat{Q_{1}} \supseteq \hat{A_{2}} \supseteq \cdots \supseteq \hat{A_{r}} \supseteq \hat{Q_{r}} \supseteq \hat{A}_{r+1}= \varnothing. $$
Assume that 
$ B=\lbrace  X_{i} \vert 0 \leq i \leq k  \rbrace$  
is an independent subset of
$ K $
satisfying
$$ \vert B \cap A_{j} \vert = \dim ( K \cap A_{j})+1 $$
for
$0 \leq j \leq r. $
If 
$ H_{i}$
is a total polar of  
$ X_{i}$
for
$ 0 \leq i \leq k$
and
$ \lbrace H_{i} \vert 0 \leq i \leq k \rbrace $
is independent, the 
$(n-k-1)$-dimensional 
subspace  
$ \displaystyle{\bigcap _{i=0}^{k} H_{i}} $
is a total polar of
$ K$
with respect to the absolute figure of the space.
\label{d4}
\end{definition}
\begin{remark}
By Definition
\ref{d4}
the subspaces 
$ P^{n}(\mathbb{R}) $
and
$\varnothing$
are total polar to each other in all Cayley-Klein spaces of dimension
$n.$
\end{remark}
The set of all total polars of a point in a Cayley-Klein plane may be a single line, a pencil of lines or all lines of the plane and the set of all total polars of a line may be a single point, a projective range of points or all points of the plane. Each of these sets is a certain subvariety of a Grassmannian called a Schubert variety.
\begin{definition}
Let
$$ \varnothing = W_{-1} \subsetneq W_{0} \subsetneq \cdots \subsetneq W_{i} \subsetneq \cdots \subsetneq W_{k} $$
be a strictly increasing sequence of subspaces in a projective space. The Schubert variety associated with this sequence is, by definition, the set of all 
$k$-dimensional 
subspaces such as
$ K $
that satisfies
$ \dim (K \cap W_{i}) \geq i $
for
$ 0\leq i \leq k .$
\end{definition}
For instance, in a projective plane
$ \textbf{P}$
the pencil of lines passing through a point
$X$
is obtained by assuming that
$ W_{0}=X $
and
$W_{1}=\textbf{P}.$
As an additional example, all lines of the plane
$ \textbf{P}$
is the Schubert variety associated to a sequence
$\varnothing \subseteq L \subseteq \textbf{P}, $
in which
$L$
is an arbitrary line of the plane. In Theorem
\ref{t2},
we show that the set of all total polars of a subspace in an arbitrary Cayley-Klein space is a Schubert variaty.
\begin{definition}
Let
$K$
be a subspace of a Cayley-Klein space with the absolute figure
$$   \hat{Q_{0}} \supseteq \hat{A_{1}} \supseteq \hat{Q_{1}} \supseteq \hat{A_{2}} \supseteq \cdots \supseteq \hat{A_{r}} \supseteq \hat{Q_{r}} \supseteq \hat{A}_{r+1}= \varnothing. $$
The polar sequence of
$K$
is a subsequence of 
$$ K^{p}=( K\cap A_{0})^{p_{0}} \supseteq \cdots \supseteq ( K\cap A_{j})^{p_{j}} \supseteq \cdots \supseteq ( K\cap A_{r})^{p_{r}} \supseteq \varnothing $$
obtained by deleting each subspace
$( K\cap A_{j})^{p_{j}}$
which is equal to
$A_{j+1}.$
\footnote{This sequence is called the system of total polar subspaces of
$K$ 
by Hermann Vogel in
\cite{vogel1994system}. 
He considers it as a generalization of the total polar of a regular subspace.
}
\end{definition}
In the next theorem, we characterize all total polars of a subspace by using its polar sequence. To prove the theorem we need some lemmas.
\begin{lemma}
Let 
$ Q $
be a quadric with vertex 
$ A $
in 
$P^{n}(\mathbb{R}).$
If  
$ H $
is a hyperplane that contains 
$ A ,$
then  
$H^{p}= X+A $
for a point
$ X $
not in
$ A. $
 \label{l1}
 \end{lemma}
\begin{lemma}
Let 
$ Q $
be a quadric with vertex
$ A $
in 
$P^{n}(\mathbb{R}).$
For a given subspace 
$ K, $
we have
\begin{enumerate}
\item
If
$ H $
is a hyperplane that contains
$ K^{p} ,$
there exists a point 
$ X \in  K \setminus A $
such that
$ H^{p}=X+A $
and so, 
$ X^{p} =H .$
\item
Let
$  H_{1}, \ldots, H_{m} $
be independent hyperplanes which are containing
$ K^{p} .$
If 
$ H_{i}^{p} = X_{i} +A  $
for a point 
$ X_{i}(\textit{\textbf{x}}_{i}) \in K \setminus A $
with
$ 1 \leq i \leq m,  $
then for any independent subset
$ Y $
of
$ A$
the set
$ Y \cup \lbrace X_{i} \vert 1 \leq i \leq m \rbrace $
is independent.
\end{enumerate}
\begin{proof} 
\begin{enumerate}
\item
The hyperplane
$ H $
contains
$ A .$
So
$ H^{p}= Z+A $
for some point
$ Z \notin A.$
Suppose that
$ Z \notin K.$
Since
$ K^{p} \subseteq H , $
$ H^{p} \subseteq (K^{p})^{p} = K+A.$
Hence,
$ Z \in X+W $
for some
$X \in K \setminus A $
and
$ W \in A. $
This gives
$$ H^{p}=Z+A \subseteq (X+W)+A= X+ A. $$ 
It follows that 
$ H^{p}= X+A $
and consequently,
$ X^{p}= H.$
\item
Let 
$Y= \lbrace X_{m+1}(\textit{\textbf{x}}_{m+1}), \ldots, X_{m+r}(\textit{\textbf{x}}_{m+r}) \rbrace $
and
$ \mathsf{M} $
be the matrix associated to
$ Q. $
Suppose that
$ \displaystyle{\sum_{i=1}^{m+r} \alpha_{i} \mathsf{X}_{i} }= 0 $
for given scalars 
$ \alpha_{i}^{,}s. $
So
$\displaystyle{\sum_{i=1}^{m+r} \alpha_{i} (\mathsf{M} \mathsf{X}_{i}) }=0,$
implying
$\displaystyle{\sum_{i=1}^{m} \alpha_{i} (\mathsf{M} \mathsf{X}_{i}) }=0.$
The hyperplanes  
$H_{i}^{,}s  $
are independent. This gives
$ \alpha_{i}=0 $
for
$ 1 \leq i \leq m,$
and then
$ \displaystyle{\sum_{i=m+1}^{m+r} \alpha_{i} \mathsf{X}_{i} }= 0 .$
 Since
$ Y $
is independent, we get
$ \alpha_{i}=0 $
for
$ 1 \leq i \leq m+r.$
\end{enumerate}
\end{proof}
\label{l3}
\end{lemma}
\begin{lemma}
Let 
$ S $
and
$ T $
be subspaces of 
$ P^{n}(\mathbb R) $
satisfying
$ S \subseteq T.$
If 
$\dim S = n-m-1 $ 
and
$\dim T = n-r-1 ,$ 
there exist hyperplanes
$ H_{0},$
$ \ldots,$
$ H_{m} $
such that
$ T=\displaystyle{\bigcap _{i=0}^{r} H_{i}} $
and
$ S=\displaystyle{\bigcap _{i=0}^{m} H_{i}}. $
Moreover,
$ \lbrace H_{i} \cap T  \vert r < i \leq m \rbrace $
is an independent set of hyperplanes of 
$ T.$
\begin{proof}
Suppose that
$ T= S+ < X_{1},\ldots,X_{m-r}> $
for an independent subset 
$ \lbrace X_{1},\ldots ,X_{m-r}\rbrace \subseteq T \setminus S. $ 
For each 
$i,$
$ 1 \leq i \leq m-r,$
let 
$ H_{i+r}$
be a hyperplane of
$ P^{n}(\mathbb R) $
that contains
$ S+ < X_{1},\ldots,\hat{X_{i}},\ldots, X_{m-r}> $
and does not contain
$ X_{i}.$
If 
$ T=\displaystyle{\bigcap _{i=0}^{r} H_{i}} $
where the hyperplanes
$ H_{0},$
$ \ldots,$
$ H_{r} $
are independent, then 
$ S=\displaystyle{\bigcap _{i=0}^{m} H_{i}} $
and 
$ \lbrace H_{i} \cap T  \vert r < i \leq m \rbrace $
is independent.
\end{proof}
\label{l4}
\end{lemma}
\begin{theorem}
In a Cayley-Klein space of dimension
$n$
with the absolute figure
$$   \hat{Q_{0}} \supseteq \hat{A_{1}} \supseteq \hat{Q_{1}} \supseteq \hat{A_{2}} \supseteq \cdots \supseteq \hat{A_{r}} \supseteq \hat{Q_{r}} \supseteq \hat{A}_{r+1}= \varnothing,$$
let
$$( K\cap A_{j_{0}})^{p_{j_{0}}} \supseteq \cdots \supseteq  ( K\cap A_{j_{i}})^{p_{j_{i}}} \supseteq \cdots \supseteq  ( K\cap A_{j_{s}})^{p_{j_{s}}}  \supseteq \varnothing$$
be the polar sequence of a
$k$-dimensional 
subspace
$K$
of
$P^{n}(\mathbb{R}).$
An
$(n-k-1)$-dimensional
subspace
$Y$
is a total polar of
$K$
if and only if for each
$i,$
$ 0 \leq i \leq s,$
$\dim (Y \cap ( K\cap A_{j_{i}})^{p_{j_{i}}} ) \geq m_{j_{i}}- k_{j_{i}}-1 $
where 
$m_{j}=\dim A_{j}$
and
$k_{j}=\dim ( K \cap A_{j})$
for
$ 0 \leq j \leq r.$ 
\begin{proof}
We first show that each total polar of
$K$
satisfies the above inequalities. Let
$ K^{\perp}=\displaystyle{\bigcap _{t=0}^{k}  H_{t}} $
be a total polar of
$K,$
where the hyperplanes
$ H_{t}$
for
$ 0 \leq t \leq k$
are total polars of points of an independent subset 
$B$
of
$ K $ 
satisfying
$ \vert B \cap A_{j} \vert = k_{j}+1 $ 
for
$0 \leq j \leq r. $
Moreover, we may assume that the hyperplanes
$H_{t}$
for
$ 0 \leq t \leq k_{j}$
are total polars of points of
$B \cap A_{j}.$
We have
$$\dim (\displaystyle{\bigcap _{t=0}^{k}  H_{t}} \cap (K\cap A_{j_{i}})^{p_{j_{i}}})
 = \dim (\displaystyle{\bigcap _{t=0}^{k_{j_{i}+1}}  H_{t}} \cap (K\cap A_{j_{i}})^{p_{j_{i}}}) \qquad
 $$
 \begin{align*}
&=\dim(\displaystyle{\bigcap _{t=0}^{k_{j_{i}+1}}  H_{t}})+ \dim(K\cap A_{j_{i}})^{p_{j_{i}}}-\dim ( (\displaystyle{\bigcap _{t=0}^{k_{j_{i}+1}}  H_{t}}) + (K\cap A_{j_{i}})^{p_{j_{i}}})\\ 
&=n-(k_{j_{i}+1}+1)+m_{j_{i}}-k_{j_{i}}+k_{j_{i}+1} -\dim ( (\displaystyle{\bigcap _{t=0}^{k_{j_{i}+1}}  H_{t}}) + (K\cap A_{j_{i}})^{p_{j_{i}}})\\
&\geq m_{j_{i}} - k_{ j_{i}}-1 
\end{align*}
\par
Now, let 
$Y$
be an
$(n-k-1)$-dimensional subspace satisfying
$\dim (Y \cap ( K\cap A_{j_{i}})^{p_{j_{i}}} ) \geq m_{j_{i}}- k_{j_{i}}-1 $
for
$ 0 \leq i \leq s.$
We show that
$Y$
is a total polar of
$K.$
Since the polar sequence of
$\varnothing$
is the sequence of vertices of the absolute figure of the space and
$\dim (P^{n}(\mathbb{R}) \cap A_{j})=m_{j}$
for 
$0 \leq j \leq r,$
the theorem holds for
$K=\varnothing. $
If
$K \neq \varnothing, $
it is proved that there exists an independent subset 
$ B=\lbrace X_{i} \vert 0 \leq i \leq k \rbrace  $
of
$ K $
that satisfies
$ \vert B \cap A_{j} \vert = k_{j}+1 $
for
$0 \leq j \leq r $
such that for each i,
$ 0 \leq i \leq k, $
there exists a total polar
$ H_{i} $
for 
$ X_{i} $
satisfying
$ Y= \displaystyle{\bigcap _{i=0}^{k}  H_{i}}.$ 
Assume that
$$ K \cap A_{l} \neq \varnothing,  \qquad K \cap A_{l+1} = \varnothing,  $$
$$ K \subseteq A_{d} \qquad \text{and}  \qquad K \nsubseteq  A_{d+1}$$
for some integers
$ 0 \leq d \leq l \leq r.$ 
Notice that
$( K \cap A_{j})^{p_{j}}= A_{j} \neq A_{j+1} $
for
$0 \leq j < d $
and
$l< j \leq r.$
It is easily seen that the polar sequence of 
$K$
includes
$( K \cap A_{l})^{p_{l}}$
and
$  ( K \cap A_{l})^{p_{l}} \subseteq Y.  $
Consider the following sequence of subspaces 
$$ Y = Y+ ( K \cap A_{l})^{p_{l}} \subseteq \cdots \subseteq Y +( K \cap A_{j})^{p_{j}} \subseteq Y+A_{j} \subseteq
\cdots \subseteq  Y+A_{d}.
$$
Applying Lemma 
\ref{l4}
for any two subsequent subspaces of the sequence gives
$$ Y= (\displaystyle{\bigcap _{i=0}^{r_{l}} H_{i})} \cap (\displaystyle{\bigcap _{i=r_{l}+1}^{n_{l-1}} H_{i})} \cap \ldots  \cap (\displaystyle{\bigcap _{i=n_{d}+1}^{r_{d}} H_{i})} \cap  (\displaystyle{\bigcap _{i=r_{d}+1}^{k} H_{i})}  $$
for hyperplanes 
$ H_{0}, \ldots, H_{k}$
where
$(i)$
for 
each
$j,$
$ d \leq j \leq l,$
$ Y+A_{j} = \displaystyle{\bigcap _{i=r_{j}+1}^{k} H_{i}}$
and
$ Y+( K \cap A_{j} )^{ p_{j}} = \displaystyle{\bigcap _{i=n_{j}+1}^{k} H_{i}} $
by assuming
$ n_{l}= -1;$
and
$(ii)$
$ \lbrace H_{i}\cap (Y+A_{j}) \mid n_{j} < i \leq r_{j} \rbrace$
for
$ d \leq j \leq l $
and
$ \lbrace H_{i}\cap (Y+(K \cap A_{j})^{p_{j}}) \mid r_{j+1} <i \leq n_{j} \rbrace$
for
$ d \leq j < l $
are independent sets of hyperplanes of 
$Y+A_{j}$
and 
$ Y+(K \cap A_{j})^{p_{j}} $ 
respectively. 
\\
Fix some integer
$j,$
$ d \leq j \leq l. $
We show that
$ \lbrace H_{i}\cap A_{j} \mid n_{j} < i \leq r_{j} \rbrace$
is an independent set of hyperplanes of
$ A_{j}.$
For each
$i,$
$ n_{j} < i \leq r_{j}, $
since
$ Y \subseteq H_{i} $
and
$ Y + A_{j} \nsubseteq H_{i}, $
$ A_{j} \nsubseteq H_{i}.$
It is easy to see that
 $(H_{i}\cap A_{j})+Y= H_{i} \cap (Y+A_{j})$
and
$$ (\displaystyle{\bigcap _{\begin{array}{c}
{ \scriptstyle t=n_{j}+1 }\\ 
{\scriptstyle t \neq i}
\end{array}}^{r_{j}}(H_{t} \cap A_{j})})+Y= \displaystyle{\bigcap _{\begin{array}{c}
{ \scriptstyle t=n_{j}+1} \\ 
{ \scriptstyle t \neq i }
\end{array}}^{r_{j}}( H_{t} \cap(Y+ A_{j}))} $$
for
$ n_{j} < i \leq r_{j}. $
From these equalities, it is followed that 
$ \lbrace H_{i}\cap A_{j} \mid n_{j} < i \leq r_{j} \rbrace$
is independent. By a similar argument, it can be shown that
$ \lbrace H_{i}\cap (K \cap A_{j})^{p_{j}}) \mid r_{j+1} <i \leq n_{j} \rbrace $
is an independent set of hyperplanes of
$ (K \cap A_{j})^{p_{j}}$
for
$ d \leq j < l. $
From the independency of these sets, it is concluded that
$ \lbrace  H_{t} \cap (K \cap A_{j_{i}})^{p_{j_{i}}}\vert 0 \leq t \leq n_{j_{i}} \rbrace $
is also an independent set of hyperplanes of 
$ (K \cap A_{j_{i}})^{p_{j_{i}}} $
for
$ d \leq i < s-(r-l-1) $
provided that
$ j_{i} \neq l.$
From
$ \dim (Y \cap  (K \cap A_{j_{i}})^{p_{j_{i}}}) \geq  m_{j_{i}}-k_{j_{i}}-1,$
it is followed that
$$ \dim ((\displaystyle{\bigcap _{j=0}^{n_{j_{i}}} H_{j}}) \cap (\displaystyle{\bigcap _{j=n_{j_{i}}+1}^{k} H_{j}}) \cap (K \cap A_{j_{i}})^{p_{j_{i}}}) \geq m_{j_{i}}-k_{j_{i}}-1.$$  
Hence
$  (m_{j_{i}}-k_{j_{i}}+k_{j_{i}+1}) - (n_{j_{i}}+1) \geq m_{j_{i}}-k_{j_{i}}-1,$
implying
$  k_{j_{i}+1} \geq n_{j_{i}} .$
Now, let
$ B_{l+1}= \varnothing$
and fix some integer
$ j$
with
$ d \leq j \leq l. $
Suppose that
$ B_{j+1}$
is an independent subset of
$ K \cap A_{j+1}$
satisfying
$ \vert B_{j+1} \cap A_{t} \vert = k_{t}+1 $
for
$ j+1 \leq t \leq l+1. $
$ \lbrace H_{i}\cap A_{j} \mid n_{j} < i \leq r_{j}\rbrace$
is an independent set of hyperplanes of
$ A_{j}$
that each of them contains
$(K \cap A_{j})^{p_{j}}. $ 
Applying Lemma 
\ref{l3}
gives an independent subset
$\lbrace X_{k_{j+1}+1},\ldots, X_{k_{j+1}+r_{j}-n_{j}} \rbrace $
of
$( K\cap A_{j} )\setminus A_{j+1}$
satisfying
$ X^{p_{j}}_{k_{j+1}+i-n_{j}}= H_{i}\cap A_{j}$
for
$ n_{j} < i \leq r_{j}. $
We consider 
$ H_{i},$
$ n_{j} < i \leq r_{j}, $
as a total polar for
$ X_{k_{j+1}+i-n_{j}}.$
The set 
$ B_{j+1} \cup \lbrace X_{k_{j+1}+1},\ldots, X_{k_{j+1}+r_{j}-n_{j}} \rbrace  $
is an independent subset of
$ K \cap A_{j}.$
Extend this set to a maximal independent subset 
$ B_{j}=\lbrace X_{0},\ldots, X_{k_{j}} \rbrace  $
of
$ K \cap A_{j}.$
In this way, we get an independent subset
$ B_{d} $
of
$ K $
satisfying
$ \vert B_{d} \cap A_{t} \vert = k_{t}+1 $
for
$ 0 \leq t \leq r $
such that we have determined individual total polars for some of its elements. For each 
$j,$
$ d< j \leq l,$
if
$ r_{j} \neq n_{j-1}, $
$(K \cap A_{j-1})^{p_{j-1}} \neq A_{j}. $
So
$ j-1 =j_{i}$
for some 
$i,$
$d \leq i < s-(r-l-1).$
Regarding the inequality 
$ n_{j_{i}} \leq k_{j_{i}+1,}$
we consider the hyperplanes
$ H_{t}$
for 
$ r_{j} < t \leq n_{j-1},$
as total polars for those points of
$ B_{j}  $
that have not yet been determined any total polar for them. The hyperplanes
$ H_{j}$
for
$ j >r_{d}, $
are considered as total polars for the remaining points of
$ B_{d}.$
\label{t1}
\end{proof}
\end{theorem}
\begin{theorem}
In a Cayley-Klein space, the set of all total polars of a subspace  is a Schubert variety.
\begin{proof}
Consider a Cayley-Klein space of dimension
$ n $
with the absolute figure 
$$   \hat{Q_{0}} \supseteq \hat{A_{1}} \supseteq \hat{Q_{1}} \supseteq \hat{A_{2}} \supseteq \cdots \supseteq \hat{A_{r}} \supseteq \hat{Q_{r}} \supseteq \hat{A}_{r+1}= \varnothing $$
and let 
$ K $
be a 
$ k$-dimensional
subspace of
$P^{n}(\mathbb{R}).$
We show that all total polars of 
$ K $
is a Schubert variety. Set
$ m_{j}=\dim A_{j}$
and
$k_{j}=\dim (K \cap A_{j})$
for
$ 0 \leq j \leq r .$
We know that
$ ( K\cap A_{j})^{p_{j}} \supseteq A_{j+1}.$
This gives
$ \dim ( K\cap A_{j})^{p_{j}} \geq \dim A_{j+1}$
and then
$  m_{j}-k_{j} -1 \geq  m_{j+1}-k_{j+1}-1$
consequently. Moreover, two sides are equal if and only if
$ ( K\cap A_{j})^{p_{j}} = A_{j+1}. $
Let
$$   ( K\cap A_{j_{0}})^{p_{j_{0}}} \supseteq \cdots \supseteq  ( K\cap A_{j_{i}})^{p_{j_{i}}} \supseteq \cdots \supseteq  ( K\cap A_{j_{s}})^{p_{j_{s}}}\supseteq \varnothing $$
be the polar sequence of
$K.$
From
$  m_{j_{i}+1}-k_{j_{i}+1}-1 =m_{j_{i+1}}-k_{j_{i+1}}-1,$
it is followed that
$$  m_{j_{0}}-k_{j_{0}}-1 > \cdots > m_{j_{i}}-k_{j_{i}}-1 > \cdots >  m_{j_{s}}-k_{j_{s}}-1 > m_{j_{s+1}}-k_{j_{s+1}}-1 = -1  . $$
It is easily seen that 
 $ m_{j_{0}}-k_{j_{0}} -1 = n-k-1. $
We set
$ W_{-1}= \varnothing$
and
$$W_{m_{j_{i}}-k_{j_{i}}-1}= ( K\cap A_{j_{i}})^{p_{j_{i}}}$$ for
$ 0 \leq i  \leq s.$ 
Regarding the equalities
$$ \dim  ( K\cap A_{j_{i}})^{p_{j_{i}}} = m_{j_{i}}-k_{j_{i}}+k_{j_{i}+1} \quad  \text{and}  \quad \dim A_{j_{i}+1}= m_{j_{i}+1},$$
for each
$ c,$ 
$$ 1 \leq  c < (m_{j_{i}}-k_{j_{i}}-1) -  (m_{j_{i+1}}-k_{j_{i+1}}-1)=(m_{j_{i}}-k_{j_{i}}+k_{j_{i}+1})- m_{j_{i}+1}, $$
assume that
$ A_{j_{i}+1} +Z_{i}^{c} $
is an
$ (m_{j_{i}+1} +c) $-dimensional
subspace of
$ ( K\cap A_{j_{i}})^{p_{j_{i}}} $ 
such that
$$  A_{j_{i}+1} +Z_{i}^{(m_{j_{i}}-k_{j_{i}}+k_{j_{i}+1})- m_{j_{i}+1}-1} \supseteq \cdots \supseteq A_{j_{i}+1} +Z_{i}^{c} \supseteq \cdots \supseteq A_{j_{i}+1} +Z_{i}^{1} . $$
For each 
$ c,$
set
$$W_{m_{j_{i+1}}-k_{j_{i+1}}-1+c}= A_{j_{i}+1} +Z_{i}^{c}.$$
We claim that the Schubert variety associated with  
$ \lbrace W_{i} \rbrace_{i=-1}^{n-k-1} $
is the Schubert variety claimed in the theorem. Due to Theorem 
\ref{t1},
it suffices to show that for a given total polar
$K^{\perp}$
of
$K$
we have
$$ \dim (K^{\perp} \cap (A_{j_{i}+1} +Z_{i}^{c})) \geq m_{j_{i+1}}-k_{j_{i+1}}-1+c.$$
From
$\dim (K^{\perp} \cap ( K\cap A_{j_{i}})^{p_{j_{i}}}) \geq m_{j_{i}}- k_{j_{i}}-1,$
it is followed that 
$\dim (K^{\perp}+ ( K\cap A_{j_{i}})^{p_{j_{i}}}) \leq n-k+k_{j_{i}+1} .$
This gives
$ \dim (K^{\perp} + A_{j_{i}+1} +Z_{i}^{c}) \leq n-k+k_{j_{i}+1}. $
Since
$ m_{j_{i}+1}-k_{j_{i}+1}=m_{j_{i+1}}-k_{j_{i+1}},$
we get
$ \dim (K^{\perp} + A_{j_{i}+1} +Z_{i}^{c}) \leq n-k+ m_{j_{i}+1}-m_{j_{i+1}}+k_{j_{i+1}}  $
and then
$ \dim (K^{\perp} \cap (A_{j_{i}+1} +Z_{i}^{c})) \geq m_{j_{i+1}}-k_{j_{i+1}}-1+c. $
\end{proof}
\label{t2}
\end {theorem}
\begin{definition}
For a given subspace 
$ K $
of a Cayley-Klein space, the set of all total polars of
$ K $
is called the polar variety of 
$ K.$
\end{definition}
Schubert varieties are non-empty. Therefore, every subspace of a Cayley-Klein space has at least one total polar. The next proposition presents sufficient and necessary condition for the uniqueness of such a total polar.
\begin{proposition}
Consider a Cayley-Klein space of dimension
$ n $
with the absolute figure 
$$   \hat{Q_{0}} \supseteq \hat{A_{1}} \supseteq \hat{Q_{1}} \supseteq \hat{A_{2}} \supseteq \cdots \supseteq \hat{A_{r}} \supseteq \hat{Q_{r}} \supseteq \hat{A}_{r+1}= \varnothing .$$
If 
 $ K $
is a subspace of  
$ P^{n}(\mathbb{R}) $ 
satisfying   
 $ K \cap A_{l} \neq \varnothing$
and 
$ K \cap A_{l+1} = \varnothing $
for an integer 
$ l $
with 
$ 0 \leq l\leq r,$ 
$K$ 
has a unique total polar 
$(K  \cap A_{l})^{p_{l}}$
if and only if it is a regular subspace.
\begin{proof}
Assume that
$$( K\cap A_{j_{0}})^{p_{j_{0}}} \supseteq \cdots \supseteq  ( K\cap A_{j_{i}})^{p_{j_{i}}} \supseteq \cdots \supseteq  ( K\cap A_{j_{s}})^{p_{j_{s}}} \supseteq \varnothing $$
is the polar sequence of 
$K$
and
$ \dim K=k. $
Let
$m_{j}=\dim A_{j}$
and 
$k_{j}=\dim (K \cap A_{j})$
for
$0 \leq j \leq r.$
Due to Theorem
\ref{t1},
it is followed that 
$ K $
has exactly one total polar if and only if
$\dim (K\cap A_{j_{0}})^{p_{j_{0}}}= n-k-1.$
In this case, the total polar of 
$K$
is 
$(K \cap A_{j_{0}})^{p_{j_{0}}}.$
Since 
$ m_{j_{0}}-k_{j_{0}}=n-k,$
it is followed that
$k_{j_{0}+1}=-1.$
This implies that
$j_{0}=l$
or
$j_{0}=l+1.$
If
$j_{0}=l,$
then
$\dim(K+A_{l})=k+m_{l}-k_{l}=n.$
Otherwise, we have 
$(K \cap A_{l})^{p_{l}}=A_{l+1}=(K \cap A_{j_{0}})^{p_{j_{0}}}.$
This gives 
$m_{l}-k_{l}=m_{l+1}-k_{l+1}=n-k,$
implying
$ \dim (K+A_{l})=n. $
In both cases, the total polar of 
$K$
equals
$(K \cap A_{l})^{p_{l}}.$
\end{proof}
\label{P1}
\end{proposition}
\begin{proposition}
Let
$ K $
and 
$ K^{\prime} $
be subspaces of a Cayley-Klein space such that
$ K \subseteq  K^{\prime}. $
Then
\begin{enumerate}
\item
Every total polar of
$ K $
contains a total polar of
$ K^{\prime}. $
\item
Every total polar of
$ K^{\prime} $
can be extended to a total polar of 
$ K. $
\end{enumerate}
\begin{proof}
Let 
$$   \hat{Q_{0}} \supseteq \hat{A_{1}} \supseteq \hat{Q_{1}} \supseteq \hat{A_{2}} \supseteq \cdots \supseteq \hat{A_{r}} \supseteq \hat{Q_{r}} \supseteq \hat{A}_{r+1}= \varnothing .$$
be the absolute figure of the space. Set
$m_{j} = \dim A_{j}$
and
$k_{j} = \dim (K \cap A_{j})$
for
$0 \leq j \leq r.$
To prove the proposition, it suffices to consider the case
$ \dim K^{\prime}-\dim K=1. $
In this case, there exists a point
$X \in (K^{\prime}\setminus K)$
and some integer 
$l,$
$0 \leq l \leq r,$
such that 
$K \cap A_{j}= K^{\prime} \cap A_{j}$
for 
$l< j \leq r+1$
and
$ (K \cap A_{j})+X= K^{\prime} \cap A_{j}$
for
$0 \leq j \leq l.$
This gives 
$\dim(K \cap A_{l})^{p_{l}}- \dim(K^{\prime} \cap A_{l})^{p_{l}}=1 $
and
$ (K \cap A_{j})^{p_{j}}= (K^{\prime} \cap A_{j})^{p_{j}} $
for 
$j \neq l.$
Let
$$( K\cap A_{j_{0}})^{p_{j_{0}}} \supseteq \cdots \supseteq  ( K\cap A_{j_{i}})^{p_{j_{i}}}=(K \cap A_{l})^{p_{l}}  \supseteq \cdots \supseteq  ( K\cap A_{j_{s}})^{p_{j_{s}}}\supseteq \varnothing $$
be the polar sequence of 
$K.$
We have
\begin{enumerate}
\item
For 
$ K^{\perp}, $
a total polar of
$ K ,$ 
extend an arbitrary maximal independent subset of
$ K^{\perp} \cap (K^{\prime} \cap A_{l})^{p_{l}} $
to a maximal independent subset of 
$ K^{\perp} \cap (K \cap A_{l})^{p_{l}}$ 
and then to a maximal independent subset of 
$ K^{\perp} \cap (K \cap A_{j_{i-1}})^{p_{j_{i-1}}}.$
Continue in this way until getting a maximal independent subset
$B $
of 
$K^{\perp} \cap ( K\cap A_{j_{0}})^{p_{j_{0}}}  =K^{\perp}.$
Suppose that
$ B^{\prime} $
is obtained by deleting one of the elements of
$ B $
added in the last step. Let
$(K^{\prime})^{\perp}$
be the subspace spanned by
$ B^{\prime}.$
This gives
$ \dim ((K^{\prime})^{\perp} \cap ( K^{\prime}\cap A_{j_{t}})^{p_{j_{t}}} ) \geq m_{j_{t}}-k_{j_{t}}-1 $
for
$i < t \leq s $
and
$\dim ((K^{\prime})^{\perp} \cap ( K^{\prime}\cap A_{j_{t}})^{p_{j_{t}}}) \geq m_{j_{t}}-k_{j_{t}}-2 $
for
$0 \leq t \leq i.$
\item
Let 
$ (K^{\prime})^{\perp} $
be a total polar of
$ K^{\prime} $
and
$r,$
$0 \leq r \leq s,$
be the maximum number that
$  (K \cap A_{j_{r}})^{p_{j_{r}}} \nsubseteq (K^{\prime})^{\perp}.$
Set
$K^{\perp}=(K^{\prime})^{\perp}+X$
for a point 
$X \in (K \cap A_{j_{r}})^{p_{j_{r}}} \setminus (K^{\prime})^{\perp}.$
It is easy to check that
$ \dim (K^{\perp}\cap ( K\cap A_{j_{t}})^{p_{j_{t}}} ) \geq m_{j_{t}}-k_{j_{t}}-1  $
for
$ 0 \leq t \leq s.  $
\end{enumerate}
\end{proof}
\label{p3}
\end{proposition}
\begin{theorem}
In a Cayley-Klein space of dimension
$ n, $
let 
$ K $
be a subspace of 
$ P^{n}(\mathbb{R}). $
If 
$ K^{\perp} $
is a total polar of 
$ K, $
then 
$ K $
is also a total polar of
$ K^{\perp}.$
\begin{proof}
 Let
$$   \hat{Q_{0}} \supseteq \hat{A_{1}} \supseteq \hat{Q_{1}} \supseteq \hat{A_{2}} \supseteq \cdots \supseteq \hat{A_{r}} \supseteq \hat{Q_{r}} \supseteq \hat{A}_{r+1}= \varnothing $$
be the absolute figure of the space. We prove the theorem by induction on dimension
$ K. $
Since the subspaces  
$ P^{n}(\mathbb{R}) $
and
$\varnothing $
are the only total polars of each other, the theorem holds for
$K=\varnothing.$
Now suppose that it is true for every 
$ (k-1) $-dimensional
subspace of
$ P^{n}(\mathbb{R}). $
Let 
$ K $
be a 
$ k $-dimensional 
subspace satisfying
$ K \cap A_{l}\neq\varnothing$
and
$ K \cap A_{l+1}=\varnothing $
for some 
$ l,$
$ 0\leq l\leq r. $
Assume that
$  K^{\perp} $
is an arbitrary total polar of
$ K. $
We show that
$ K$
is a total polar of
$  K^{\perp} .$
If
$K^{\perp} \nsupseteq A_{l},$
then 
$K \cap (K^{\perp} \cap A_{l})^{p_{l}}\neq \varnothing .$
Let
$ X $
be a point of
$K \cap (K^{\perp} \cap A_{l})^{p_{l}}$
in this case and a point of
$ K \cap A_{l} $
otherwise. Assume that
$ K-X $
is a complement of 
$X$
in 
$K.$
Due to Proposition
\ref{p3}
there exists a total polar of
$ K-X $
that contains
$ K^{\perp}. $
Let 
$  K^{\perp}+Z $
be such a total polar for some point
$ Z \notin K^{\perp}. $
By the induction hypothesis, we conclude that
$ K-X $
is also a total polar of
$  K^{\perp}+Z. $
Repeated application of Proposition
\ref{p3}
ensures the existence of a total polar 
$ K^{\prime} $
of
$  K^{\perp} $
that contains
$ K-X. $
After this preliminary step, we return to prove the assertion. We know that
$ K^{\perp} \supseteq (K \cap A_{l})^{p_{l}} \supseteq A_{l+1}. $
Let
$ K^{\perp} \supseteq  A_{j+1} $
and
$ K^{\perp} \nsupseteq A_{j}$
for some 
$ j \leq l. $
It is easily seen that 
$  (K^{\perp} \cap A_{j})^{p_{j}}\neq A_{j+1} $
and
$ (K^{\perp} \cap A_{i})^{p_{i}}= A_{i+1}$
for
$ j< i.$
Let
$$( K^{\perp}\cap A_{j_{0}})^{p_{j_{0}}} \supseteq \cdots  \supseteq  ( K^{\perp} \cap A_{j_{s}})^{p_{j_{s}}}= (K^{\perp} \cap A_{j})^{p_{j}} \supseteq \varnothing $$
be the polar sequence of
$K^{\perp}$
and
$d_{j}=\dim (K^{\perp} \cap A_{j})$
for
$ 0 \leq j \leq r. $
Since
$ K^{\prime} $
is a total polar of
$ K^{\perp} $
and contains
$ K-X ,$
it is followed that
$\dim(K ^{\prime}\cap ( K^{\perp} \cap A_{j_{i}})^{p_{j_{i}}}) \geq m_{j_{i}}-d_{j_{i}}-1 $
and consequently,
$\dim((K-X)\cap ( K^{\perp} \cap A_{j_{i}})^{p_{j_{i}}}) \geq m_{j_{i}}-d_{j_{i}}-2 $
for
$ 0 \leq i \leq s. $
On the other hand
$ X \in K \cap ( K^{\perp} \cap A_{j_{i}})^{p_{j_{i}}}$
for
$ 0 \leq i \leq s,$
implying that
$\dim(K \cap ( K^{\perp} \cap A_{j_{i}})^{p_{j_{i}}}) \geq m_{j_{i}}-d_{j_{i}}-1 $
for
 $ 0 \leq i \leq s.$
\end{proof}
\end{theorem}
In the next proposition all total polars of a hyperplane of a Cayley-Klein space are determined. 
\begin{proposition}
Let 
$ H $
be a hyperplane of a Cayley-Klein space with the absolute figure
$$   \hat{Q_{0}} \supseteq \hat{A_{1}} \supseteq \hat{Q_{1}} \supseteq \hat{A_{2}} \supseteq \cdots \supseteq \hat{A_{r}} \supseteq \hat{Q_{r}} \supseteq \hat{A}_{r+1}= \varnothing .$$
Suppose
$ H $
contains 
$ A_{i} $
and does not contain   
$ A_{i-1} $
for some 
$ i, \ 0< i \leq r+1. $
A point
 $ X $
is a total polar of
$ H$
if and only if
$ X \in (H \cap A_{i-1})^{p_{i-1}}. $
\begin{proof}
We first show that
$(H \cap A_{j})^{p_{j}}= A_{j+1} $
for 
$ j\neq i-1. $
Set
$ m_{j}=\dim A_{j} $
for
$ 0 \leq j \leq r+1. $
If 
$ j< i-1 ,$
then
$$ \dim (H \cap A_{j})^{p_{j}} = m_{j} - \dim (H \cap A_{j}) + \dim (H \cap A_{j+1}).$$
Since
$ H $
does not contain neither
$ A_{j} $
nor
$ A_{j+1}, $
we get
$$ \dim (H \cap A_{j})^{p_{j}}= m_{j}-(m_{j}-1)+m_{j+1}-1=m_{j+1} .$$
$(H \cap A_{j})^{p_{j}}$
contains 
$ A_{j+1}$
and so, equals it. If  
$ i-1< j, $
then 
$ H $
contains
$ A_{j} $
and it is obvious that
$(H \cap A_{j})^{p_{j}}=A_{j+1}.$
For 
$ j=i-1, $
it is concluded that
$(H \cap A_{j})^{p_{j}} \neq A_{j+1}$
due to Lemma
\ref{l1}.
In view of Theorem 
\ref{t1}
and above considerations, the validness of the proposition can be easily verified.
\end{proof}
\label{p4}
\end{proposition}
In the remaining part of this section, we compare our definition for total polars with the other definitions introduced in section two. It is easily seen that in each Cayley-Klein plane, the set of all total polars of a point or a line of the plane equals to the set of polars corresponded to it in
\cite{richter2011perspectives}.
We show that the set of all total polars of a subspace in a Cayley-Klein space includes all polars associated to it in
\cite{struve2004projective}.
Consider a Cayley-Klein space of dimension 
$ n $
with the absolute figure 
$$   \hat{Q_{0}} \supseteq \hat{A_{1}} \supseteq \hat{Q_{1}} \supseteq \hat{A_{2}} \supseteq \cdots \supseteq \hat{A_{r}} \supseteq \hat{Q_{r}} \supseteq \hat{A}_{r+1}= \varnothing. $$
Let
$K$
be a 
$k$-dimensional 
subspace of 
$ P^{n}(\mathbb{R})$
satisfying
$K \cap A_{l} \neq \varnothing $
and
$K \cap A_{l+1}= \varnothing $
for some integer
$l,$
$0 \leq l \leq r.$
If
$K^{\perp}$ 
is a polar of
$K$
according to
\cite{struve2004projective},
then
$ ( K \cap A_{j})^{p_{j}} = (K^{\perp} \cap A_{j})+A_{j+1}$
for
$0 \leq j \leq r.$
Since
$ K \cap A_{j}=\varnothing $
for
$ l < j \leq r,$
it is followed that 
$ K^{\perp} \supseteq ( K \cap A_{l})^{p_{l}} .$
This implies that
$ \dim (K^{\perp} \cap ( K \cap A_{l})^{p_{l}}) = m_{l}-k_{l}-1.$
Now, assume that 
$\dim (K^{\perp} \cap ( K \cap A_{j})^{p_{j}}) \geq m_{j}-k_{j}-1$
for an integer
$j,$
$0 < j \leq l.$
We have
\begin{align*}
 & \dim (K^{\perp} \cap ( K \cap A_{j-1})^{p_{j-1}}) 
 \\
 & = \dim  K^{\perp}+ \dim ( K \cap A_{j-1})^{p_{j-1}}- \dim( K^{\perp}+  (K \cap A_{j-1})^{p_{j-1}}) \\
 & = \dim  K^{\perp}+ \dim ( K \cap A_{j-1})^{p_{j-1}} - \dim ( K^{\perp}+ (K^{\perp} \cap A_{j-1})+A_{j}) \\
 &= \dim ( K \cap A_{j-1})^{p_{j-1}} - \dim A_{j}+\dim ( K^{\perp} \cap A_{j})\\
 & \geq m_{j-1}-k_{j-1}+k_{j}-m_{j}+m_{j}-k_{j}-1 =m_{j-1}-k_{j-1}-1
\end{align*}
Regarding Theorem
\ref{t1},
it is concluded that 
$K^{\perp}$
is a total polar of 
$K.$
\section{Tangent Cones and Spheres}
In this section, we define tangent subspaces to the absolute figure of a Cayley-Klein space. We will determine all tangent lines through a given point not in the second vertex. It is proved that the set of all points on these lines is a quadric called a tangent cone. Finally, we obtain the equations of all spheres using tangent cones.
\begin{definition}
Consider a Cayley-Klein space of dimension 
$ n .$
A subspace 
$ K $
of 
$ P^{n}(\mathbb{R})$
is tangent to the absolute figure of the space if there exists a total polar 
$ K^{\perp} $
of 
$ K $
with
$ K \cap K^{\perp} \neq \varnothing. $
Briefly, we say that 
$ K $
is a tangent subspace.
\end{definition}
\begin{proposition}
Consider a Cayley-Klein space with the absolute figure
$$   \hat{Q_{0}} \supseteq \hat{A_{1}} \supseteq \hat{Q_{1}} \supseteq \hat{A_{2}} \supseteq \cdots \supseteq \hat{A_{r}} \supseteq \hat{Q_{r}} \supseteq \hat{A}_{r+1}= \varnothing, $$
such that
$ A_{1} $
is a hyperplane. Let
$X$
be a point not in
$ A_{1}. $
A line
$L$
through
$X$
is tangent to the absolute figure if and only if it passes through a point of
$ Q_{1}. $
\begin{proof}
The line
$ L $
intersects 
$ A_{1} $
at a point
$ Y. $
Suppose that 
$ Y \in A_{k} \setminus A_{k+1} $
for some integer 
$ k \geq 1. $
We show that 
$ L $
is a tangent line if and only if
$ Y \in Q_{1}.$
If 
$ k=1, \ L $
has a unique total polar 
$ Y^{p_{1}}. $
So
$ L $
is tangent to the absolute figure if and only if 
$ Y^{p_{1}}$
contains
$ Y .$
That is 
$ Y \in Q_{1}. $
If
$ k>1,$
then every hyperplane of 
$ A_{1}$  
containing
$ A_{2}$ 
is a total polar of
$ L$ 
by Theorem
\ref{t1}.
Such a total polar intersects 
$ L $
at 
$Y.$
\end{proof}
\label{p5}
\end{proposition}
\begin{proposition}
Consider a Cayley-Klein space with the absolute figure
$$   \hat{Q_{0}} \supseteq \hat{A_{1}} \supseteq \hat{Q_{1}} \supseteq \hat{A_{2}} \supseteq \cdots \supseteq \hat{A_{r}} \supseteq \hat{Q_{r}} \supseteq \hat{A}_{r+1}= \varnothing ,$$
such that 
$ A_{1} $
is not a hyperplane. Let 
$ X $
be a point not in
$A_{1}$
and 
$ L $
a line through it. Then
\begin{enumerate}
\item
For
$ X \notin Q_{0}, \ L $
is a tangent line if and only if it passes through a point of
$ X^{p} \cap Q_{0}.$
\item
For 
$ X \in Q_{0}, \ L $
is a tangent line if and only if
$ L \subseteq X^{p}. $
\end{enumerate}
\begin{proof}
\begin{enumerate}
\item
Since
$ X \notin Q_{0},$
$ X \notin X^{p}.$
So
$ L $
intersects 
$ X^{p}$
at a point such as
$ Y. $
Similar to the proof of the previous proposition we can show that 
$ L $
is a tangent line if and only if
$ Y \in X^{p} \cap Q_{0}. $
\item
Assume that 
$ L \nsubseteq X^{p}. $
We show that
$ L $
is not a tangent line. Since
$ A_{1} \subseteq X^{p}, $
$ L \cap A_{0} \neq \varnothing$
and
$L \cap A_{1} = \varnothing.$
Thus
$L$
is regular and has a unique total polar
$ X^{p} \cap Y^{p} $
for some point
$ Y \neq X $
of 
$ L.$
If
$ L $
is a tangent line, it intersects
$ X^{p} \cap Y^{p}$
at 
$ X. $
This gives
$ X \in Y^{p} $
and we get
$ Y \in X^{p} $
which is a contradiction. Now, Suppose that $ L \subseteq X^{p}. $
We consider two cases 
$ L \cap A_{1} = \varnothing $
and
$ L \cap A_{1} \neq \varnothing. $
If 
$ L \cap A_{1} = \varnothing ,$
$ L$
has a unique total polar
$ X^{p} \cap Y^{p} $
for some point
$ Y \neq X $
of 
$ L.$
Since
$ X \in X^{p} \cap Y^{p},$
$ L $
intersects its total polar and is tangent to the absolute figure of the space. If
$  L \cap A_{1} \neq \varnothing,  $
$L$
intersects
$ A_{1}$
at a point 
$ Y.$
Regarding Theorem
\ref{t1},
it is concluded that any hyperplane of
$ X^{p} $
containing
$ A_{1} $
is a total polar of
$ L. $
Such a total polar intersects
$L$
at
$Y.$
\end{enumerate}
\end{proof}
\end{proposition}
\begin{remark}
In a Cayley-Klein space of dimension
$ n $
with the absolute figure
$$   \hat{Q_{0}} \supseteq \hat{A_{1}} \supseteq \hat{Q_{1}} \supseteq \hat{A_{2}} \supseteq \cdots \supseteq \hat{A_{r}} \supseteq \hat{Q_{r}} \supseteq \hat{A}_{r+1}= \varnothing,$$
set 
$ n_{i}=\dim A_{i} - \dim A_{i+1} $
for
$ 0 \leq i \leq r. $
Let 
$ \mathsf{E}_{i} $
be a diagonal matrix of order
$ n_{i} $
where the diagonal entries equal 
$ 1 $
or
$ -1. $
Assume that the number of 
$ -1 $
elements is
$ q_{i}. $
There exists a projective coordinate system for 
$ P^{n}(\mathbb{R}) $
in which the quadrics of the absolute figure and their vertices have the following equations:
 \begin{align*}
Q_{0} &: \mathsf{X}^{T}_{0}\mathsf{E}_{0}\mathsf{X}_{0}=x^{2}_{0}+\cdots+
x^{2}_{n_{0}-q_{0}-1}-x^{2}_{n_{0}-q_{0}}-\cdots-x^{2}_{n_{0}-1}=0;  \\
\mathsf{X}_{0} &=(\begin{array}{c c c }
x_{0} & \ldots & x_{n_{0}-1} 
\end{array})^{T}. \\    
A_{1}  &: x_{0}=\cdots=x_{n_{0}-1}=0.\\ 
 Q_{1} & :  \mathsf{X}^{T}_{1}\mathsf{E}_{1}\mathsf{X}_{1}=x^{2}_{n_{0}}+\cdots+x^{2}_{n_{0}+n_{1}-q_{1}-1}-x^{2}_{n_{0}+n_{1}-q_{1}}-\cdots-x^{2}_{n_{0}+n_{1}-1}=0, \\
 \mathsf{E}_{0}&\mathsf{X}_{0}=\mathsf{O}_{0}, \
  \mathsf{X}_{1}=(\begin{array}{c c c }
x_{n_{0}} & \ldots & x_{n_{0}+n_{1}-1} 
\end{array})^{T}.\\
A_{2} &: x_{0}=\cdots=x_{n_{0}-1}=\cdots=x_{n_{0}+n_{1}-1}=0. \\
\vdots
\end{align*}
$$
Q_{r}: \mathsf{X}^{T}_{r}\mathsf{E}_{r}\mathsf{X}_{r}=x^{2}_{n_{0}+ \cdots +n_{r-1}}+\cdots+x^{2}_{n_{0}+ \cdots+ n_{r}-q_{r}-1} -x^{2}_{n_{0}+ \cdots + n_{r}-q_{r}} -\cdots \quad  $$
\begin{align*}
 &  - x^{2}_{n_{0} + \cdots+ n_{r}-1} =0 , \
 \mathsf{E}_{0} \mathsf{X}_{0}=\mathsf{O}_{0},\   \mathsf{E}_{1} \mathsf{X}_{1}=\mathsf{O}_{1}, \ldots ,  \mathsf{E}_{r-1} \mathsf{X}_{r-1}=\mathsf{O}_{r-1}; \quad \quad \quad \\
 & \mathsf{X}_{r}=(\begin{array}{c c c }
x_{n_{0}+ \cdots +n_{r-1}} & \ldots& x_{n}
\end{array})^{T}.
\end{align*}
This coordinate system is called the normal projective coordinate system.
\label{r2}
\end{remark}
\begin{theorem}
In a Cayley-Klein space of dimension
$ n $
with the absolute figure
$$   \hat{Q_{0}} \supseteq \hat{A_{1}} \supseteq \hat{Q_{1}} \supseteq \hat{A_{2}} \supseteq \cdots \supseteq \hat{A_{r}} \supseteq \hat{Q_{r}} \supseteq \hat{A}_{r+1}= \varnothing,$$
let 
$ Z( \textit{\textbf{z}}) $
be a point not in
$ A_{1}. $
The set of all points on tangent lines through
$ Z $
is a quadric.
\begin{proof}
Consider the normal projective coordinate system for the space. We first assume that 
$ A_{1} $
is a hyperplane. For a given point
$ Y( \textit{\textbf{y}}) \neq Z, $
the line 
$ Y+Z $
is a tangent line if and only if it passes through
$ Q_{1} .$
It means that
$ (Y+tZ)( \textit{\textbf{y}}+t  \textit{\textbf{z}}) $
is a point of 
$ Q_{1} $
for a scalar
$ t. $
For some integer
$n_{1},$
$ 1 \leq  n_{1} \leq n, \ Q_{1} $
has the equation
$ x_{0}=0, \  \mathsf{X}^{T}_{1}\mathsf{E}_{1}\mathsf{X}_{1}=0 $
where
$ \mathsf{X}_{1}=(\begin{array}{c c c }
x_{1} & \ldots & x_{n_{1}} 
\end{array})^{T}.$
Consequently
$ (\mathsf{Y}_{1}-y_{0} \mathsf{Z}_{1})^{T} \mathsf{E}_{1} (\mathsf{Y}_{1}-y_{0} \mathsf{Z}_{1})=0$
in which
 $ \mathsf{Z}_{1}=(\begin{array}{c c c }
z_{1} & \ldots & z_{n_{1}} 
\end{array})^{T}$
and
  $ \mathsf{Y}_{1}=(\begin{array}{c c c }
y_{1} & \ldots & y_{n_{1}} 
\end{array})^{T}.$
This gives
$$ y_{0}^{2} \mathsf{Z}^{T}_{1}\mathsf{E}_{1}\mathsf{Z}_{1} +  \mathsf{Y}^{T}_{1}\mathsf{E}_{1}\mathsf{Y}_{1} -2 y_{0} \mathsf{Y}^{T}_{1}\mathsf{E}_{1}\mathsf{Z}_{1}=0.  $$
It is a homogeneous equation of degree 2 describing a quadric in 
$ P^{n}(\mathbb{R}) .$
Now, suppose that 
$ A_{1} $
is not a hyperplane. Let 
$ Z \notin Q_{0}  $
and 
$ Y( \textit{\textbf{y}}) $
be a point with
$ Y\neq Z. $
The line 
$ Y+Z $
is a tangent line if and only if it passes through
$ Z^{p} \cap Q_{0}. $
It means that
$ (Y+tZ) ( \textit{\textbf{y}}+t \textit{\textbf{z}}) \in Z^{p} \cap Q_{0} $
for a scalar
$ t. $
We have
$$ Q_{0}:  \mathsf{X}^{T}_{0}\mathsf{E}_{0}\mathsf{X}_{0}=0, \quad \mathsf{X}_{0}=(\begin{array}{c c c }
x_{0} & \ldots & x_{n_{0}-1} 
\end{array})^{T}$$
for some integer
$ n_{0}. $
Since
$ Y+tZ \in Z^{p}, $
$ (\mathsf{Y}_{0}+t \mathsf{Z}_{0})^{T} \mathsf{E}_{0} \mathsf{Z}_{0}=0 $
where 
$ \mathsf{Y}_{0}=(\begin{array}{c c c }
y_{0} & \ldots & y_{n_{0}-1} 
\end{array})^{T}$
and
$ \mathsf{Z}_{0}=(\begin{array}{c c c }
z_{0} & \ldots & z_{n_{0}-1} 
\end{array})^{T}.$
It gives
$t=- \dfrac{\mathsf{Y}^{T}_{0}\mathsf{E}_{0}\mathsf{Z}_{0}}{\mathsf{Z}^{T}_{0}\mathsf{E}_{0}\mathsf{Z}_{0}}.$
From
$ Y+tZ \in Q_{0}, $
it is followed that 
$ (\mathsf{Y}_{0}+t \mathsf{Z}_{0})^{T} \mathsf{E}_{0} (\mathsf{Y}_{0}+t \mathsf{Z}_{0}) =0. $
Therefore,
$$ (\mathsf{Z}^{T}_{0}\mathsf{E}_{0}\mathsf{Z}_{0}) (\mathsf{Y}^{T}_{0}\mathsf{E}_{0}\mathsf{Y}_{0}) - (\mathsf{Y}^{T}_{0}\mathsf{E}_{0}\mathsf{Z}_{0})^{2}=0$$
that is a homogeneous equation of degree 2. If
$ Z \in Q_{0}, $
the desired point set is equal to the hyperplane 
$ Z^{p} $
which is a quadric in
$ P^{n}(\mathbb{R}) .$
\end{proof}
\label{t3}
\end{theorem}
\begin{definition}
Suppose that the assumptions of Theorem 
\ref{t3}
hold. The set of all points on tangent lines through 
$ Z $
is called the tangent cone with vertex 
$ Z. $
\end{definition}
\begin{definition}
Consider a Cayley-Klein space with the absolute figure   
$$   \hat{Q_{0}} \supseteq \hat{A_{1}} \supseteq \hat{Q_{1}} \supseteq \hat{A_{2}} \supseteq \cdots \supseteq \hat{A_{r}} \supseteq \hat{Q_{r}} \supseteq \hat{A}_{r+1}= \varnothing .$$
Let 
$Z$
be a point in 
$ P^{n}(\mathbb{R}) \setminus Q_{0}.$
Suppose that 
$ T_{Z} $
is the tangent cone with vertex 
$ Z. $
Every quadric in the pencil spanned by 
$ Q_{0} $
and 
$ T_{Z} $
is called a sphere with center 
$ Z. $
\label{d8}
\end{definition}
Under the assumptions of Definition
\ref{d8},
a quadric 
$ \lambda Q_{0} + \mu T_{Z}  $
is a sphere with center
$ Z $
for some real numbers 
$ \lambda $
and 
$ \mu. $
Suppose that 
$ A_{1} $
is not a hyperplane. In the normal projective coordinate system, 
$ Q_{0} $
and 
$ T_{Z} $
have the following equations
$$ Q_{0}:  \mathsf{X}^{T}_{0}\mathsf{E}_{0}\mathsf{X}_{0} =0, \quad T_{Z}: (\mathsf{Z}^{T}_{0}\mathsf{E}_{0}\mathsf{Z}_{0}) (\mathsf{X}^{T}_{0}\mathsf{E}_{0}\mathsf{X}_{0}) - (\mathsf{X}^{T}_{0}\mathsf{E}_{0}\mathsf{Z}_{0})^{2}=0$$
where 
$ \mathsf{X}_{0}=(\begin{array}{c c c }
x_{0} & \ldots & x_{n_{0}-1} 
\end{array})^{T} $
and 
$ \mathsf{Z}_{0}=(\begin{array}{c c c }
z_{0} & \ldots & z_{n_{0}-1} 
\end{array})^{T}. $
Therefore, the equation of 
$ \lambda Q_{0} + \mu T_{Z}  $
is as follows
$$ \lambda (\mathsf{X}^{T}_{0}\mathsf{E}_{0}\mathsf{X}_{0}) + \mu ((\mathsf{Z}^{T}_{0}\mathsf{E}_{0}\mathsf{Z}_{0}) (\mathsf{X}^{T}_{0}\mathsf{E}_{0}\mathsf{X}_{0}) - (\mathsf{X}^{T}_{0}\mathsf{E}_{0}\mathsf{Z}_{0})^{2})=0.$$
This gives
$$1+ \dfrac{ \lambda}{\mu (\mathsf{Z}^{T}_{0}\mathsf{E}_{0}\mathsf{Z}_{0})}= \dfrac{ (\mathsf{X}^{T}_{0}\mathsf{E}_{0}\mathsf{Z}_{0})^{2} }{ (\mathsf{Z}^{T}_{0}\mathsf{E}_{0}\mathsf{Z}_{0}) (\mathsf{X}^{T}_{0}\mathsf{E}_{0}\mathsf{X}_{0})}. $$
Suppose that
$ 1+ \dfrac{ \lambda}{\mu (\mathsf{Z}^{T}_{0}\mathsf{E}_{0}\mathsf{Z}_{0})} $
equals
$ \alpha^{2} $
or
$ (\alpha i)^{2} $
for some real number
$ \alpha >0. $
Regarding the distance formulas for points of the space, the radius of the sphere 
$ \lambda Q_{0} + \mu T_{Z}  $
is defined to be the complex number
$ \arccos (\alpha) $
or
$ \arccos (\alpha i). $
If 
$ A_{1} $
is a hyperplane, we have
$$ Q_{0}: x_{0}^{2}=0, \quad  T_{Z}: x_{0}^{2} (\mathsf{Z}^{T}_{1}\mathsf{E}_{1}\mathsf{Z}_{1}) +  \mathsf{X}^{T}_{1}\mathsf{E}_{1}\mathsf{X}_{1} -2 x_{0} (\mathsf{X}^{T}_{1}\mathsf{E}_{1}\mathsf{Z}_{1})=0.$$
Therefore, the sphere 
$ \lambda Q_{0} + \mu T_{Z}  $
has the following equation
$$  \dfrac{\lambda}{\mu} x_{0}^{2} + (\mathsf{Z}^{T}_{1}\mathsf{E}_{1}\mathsf{Z}_{1})  x_{0}^{2} + \mathsf{X}^{T}_{1}\mathsf{E}_{1}\mathsf{X}_{1} -2 x_{0} (\mathsf{X}^{T}_{1}\mathsf{E}_{1}\mathsf{Z}_{1})=0. $$
If
$ \dfrac{-\lambda}{\mu}$
equals
$ \alpha^{2} $
or
$ (\alpha i)^{2} $
for some real number
$ \alpha >0, $
the radius of the sphere is defined to be
$ \alpha $
or
$ \alpha i. $
In this way, we are able to obtain all spheres of a Cayley-Klein space without using the metric of the space.
\section{Reflections}
Every reflection in a Cayley-Klein space is defined by means of two subspaces which are total polar to each other. In this section, we give the formal definition of this notion and prove that every reflection of a Cayley-Klein space is a motion of the space. Also we show that in a Cayley-Klein space of dimension 
$n$
each motion is a product of at most 
$n+1$
reflections in point-hyperplane pairs.
\begin{definition}
Let
$K$
and
$K^{\prime}$
be two disjoint subspaces of
$P^{n}(\mathbb{R})$
satisfying
$\dim K + \dim K^{\prime}=n-1.$
For a point
$X$
not in
$K \cup K^{\prime},$
there exists a unique line passing through 
$X$
and intersecting
$K$
and 
$K^{\prime}.$
The intersections of the line with 
$K$
and
$K^{\prime}$
are called the projection of
$X$
onto
$K$
in the direction of
$K^{\prime}$
and the projection of
$X$
onto
$K^{\prime}$
in the direction of
$K$
respectively.
\end{definition}
\begin{definition}
Let
$K$
and
$K^{\prime}$
be two nonempty disjoint subspaces of
$P^{n}(\mathbb{R})$
satisfying
$\dim K + \dim K^{\prime}=n-1.$
The involution in the pair
$(K,K^{\prime})$
is a collineation of
$P^{n}(\mathbb{R})$
that fixes each point of
$K \cup K^{\prime}$
and maps any other point
$X$
to the harmonic conjugate of
$X$
with respect to the two projections of
$X$
onto
$K$
and 
$K^{\prime}.$ 
If 
$K$
and 
$K^{\prime}$
are total polar with respect to the absolute figure of a Cayley-Klein space, this involution is called the reflection in the pair
$(K,K^{\prime}).$
\end{definition}
Assume that
$ \lbrace  X_{0}, \ldots,  X_{k}  \rbrace$  
is an independent set of points of
$K$
and 
$ \lbrace  H_{0}, \ldots,  H_{k}  \rbrace$ 
is an independent set of hyperplanes containing
$K^{\prime}.$
Suppose that in a projective coordinate system,
$(\begin{array}{c c c }
x^{i}_{0} & \ldots & x^{i}_{n} 
\end{array})^T$
and
$(\begin{array}{c c c }
b^{0}_{j} & \ldots & b^{n}_{j} 
\end{array})$
are the coordinate matrices of
$ X_{i}$
and 
$ H_{j} $
respectively for
$ 0 \leq i, j\leq k.$
By setting 
$$
\mathsf{K}= \begin{pmatrix}
  x_{0}^{0}  &  x_{0}^{1} & \ldots &  x_{0}^{k} \\
   \vdots & \vdots &  & \vdots \\
  x_{n}^{0}  & x_{n}^{1} & \ldots & x_{n}^{k}
\end{pmatrix} \text{and} \ 
\mathsf{K^{\prime}}= \begin{pmatrix}
b^{0}_{0} & \ldots &  b^{n}_{0} \\ 
b^{0}_{1} & \ldots &  b^{n}_{1} \\
\vdots  &  & \vdots  \\ 
b^{0}_{k} & \ldots &  b^{n}_{k} 
\end{pmatrix}, 
$$
the matrix of the involution in the pair
$(K, K^{\prime})$
is equal to
$ \mathsf{I}_{n+1}-2 \mathsf{K} (\mathsf{K}^{\prime} \mathsf{K})^{-1} \mathsf{K}^{\prime}  $
\cite{yaglom1964projective}.
Consider a Cayley-Klein space of dimension
$n$
with the absolute figure
$$   \hat{Q_{0}} \supseteq \hat{A_{1}} \supseteq \hat{Q_{1}} \supseteq \hat{A_{2}} \supseteq \cdots \supseteq \hat{A_{r}} \supseteq \hat{Q_{r}} \supseteq \hat{A}_{r+1}= \varnothing. $$
The matrix associated to a motion of the space in the normal projective coordinate system, introduced in Remark 
\ref{r2},
is a block lower triangular matrix of order
$n+1$
such as 
 $$\begin{pmatrix}
 \mathsf{U}_{0} & & \bigzero \\
   & \ddots  &   \\
* &  & \mathsf{U}_{r} \\
\end{pmatrix}
$$
where
$ \mathsf{U}_{i} $
is a matrix of dimension 
$ n_{i}$
satisfying
$ \mathsf{U}_{i}^{T} \mathsf{E}_{i} \mathsf{U}_{i}= \mathsf{E}_{i}$
for
$0 \leq i \leq r.$
\begin{theorem}
Every reflection in a Cayley-Klein space is a motion of the space.
\begin{proof}
Consider a Cayley-Klein space of dimension
$n$
with the absolute figure
$$   \hat{Q_{0}} \supseteq \hat{A_{1}} \supseteq \hat{Q_{1}} \supseteq \hat{A_{2}} \supseteq \cdots \supseteq \hat{A_{r}} \supseteq \hat{Q_{r}} \supseteq \hat{A}_{r+1}= \varnothing. $$
Let
$K$
be a
$k$-dimensional
subspace of the space and
$K^{\perp}$
a total polar of
$K$
satisfying
$ K \cap K^{\perp}=\varnothing.$
We show that the reflection in the pair
$(K, K^{\perp})$
is a motion of the space. We have 
$ K \cap A_{l} \neq \varnothing $
and 
$ K \cap A_{l+1} = \varnothing $
for some integer 
$ 0 \leq l \leq r.$
Assume that
$ K \subseteq A_{d}$
and 
$ K \nsubseteq A_{d+1}$
for some 
$ 0 \leq d \leq r.$
Set 
$m_{j}=\dim A_{j}$
and
$k_{j}=\dim (K \cap A_{j})$
for
$ 0 \leq j \leq r.$
Let 
$ B= \lbrace X_{i} \vert 0 \leq i \leq k  \rbrace $ 
be an independent subset of 
$K $
satisfying 
$ B \cap A_{j} = \lbrace  X_{0}, \ldots,  X_{k_{j}}  \rbrace$
for 
$ d \leq j \leq l$
and
$ K^{\perp} = \displaystyle{\bigcap_{i=0}^{k} H_{i}},$
where 
$ H_{i}$
is a total polar of 
$ X_{i}$
for
$ 0 \leq i \leq k.$
First, we show that for each
$j,$
$d \leq j \leq l$
and each 
$i,$
$k_{j+1} < i \leq k_{j},$
the hyperplane
$H_{i}$
does not contain
$ A_{j}.$
We will prove this by induction on 
$j .$
Assume that the number of hyperplanes
$H_{i}$
with
$0 \leq i \leq k_{l}$
not containing 
$A_{l}$
equals
$k_{l}+1-t$
for some integer
$t \geq 0. $
These hyperplanes intersect 
$ A_{l}$
in
$ k_{l}+1-t $
independent hyperplanes
$ \lbrace X_{i}^{p_{l}} \vert H_{i}\nsupseteq A_{l} , \ 0 \leq i \leq k_{l} \rbrace$
of
$ A_{l}.$
So
$$ \dim (K^{\perp} \cap A_{l} ) = \dim ( \displaystyle{\bigcap _{i=0}^{k_{l}} H_{i}} \cap A_{l})= m_{l} -( k_{l}+1-t ). $$
Since 
$ K \cap K^{\perp}=\varnothing,$
it is followed that
$ (K \cap A_{l}) \cap (K^{\perp} \cap A_{l})=\varnothing.$
This gives
$$ k_{l}+ m_{l}-( k_{l}+1-t ) - \dim (( K \cap A_{l}) + (K^{\perp} \cap A_{l}))=-1, $$
implying
$  t=0.$
Now suppose that for each 
$ j,$
$ s < j  \leq l,$
we have 
$H_{i}\nsupseteq A_{j} $
for 
$ k_{j+1} < i \leq k_{j}. $
To complete the proof it suffices to show that 
$H_{i}\nsupseteq A_{s} $
for
$ k_{s+1} < i \leq k_{s}. $
Notice that for each 
$j,$
$ s < j  \leq l,$
the set
$$ \lbrace H_{i} \cap A_{j} \vert k_{j+1} < i \leq k_{j}\rbrace =  \lbrace X_{i}^{p_{j}} \vert k_{j+1} < i \leq k_{j}\rbrace$$
is an independent set of hyperplanes of 
$ A_{j}$
which are containing 
$ A_{j+1}.$
This implies that the set 
$  \lbrace H_{i} \cap A_{s} \vert H_{i}\nsupseteq A_{s}  , \ 0 \leq i \leq k_{s} \rbrace, $ 
which is equal to
$$  \lbrace X_{i}^{p_{s}} \vert H_{i}\nsupseteq A_{s} , k_{s+1} < i \leq k_{s} \rbrace \cup \lbrace H_{i} \cap A_{s} \vert 0 \leq i \leq k_{s+1} \rbrace,  $$
is an independent set of hyperplanes of
$ A_{s}.$
From 
$ (K \cap A_{s}) \cap (K^{\perp} \cap A_{s})=\varnothing,$
it follows that 
$ (K \cap A_{s}) \cap (\displaystyle{\bigcap _{i=0}^{k_{s}} H_{i}} \cap A_{s})=\varnothing.$
If
$t$
is the number of hyperplanes 
$H_{i}$
not containing 
$A_{s}$
for
$ k_{s+1} < i \leq k_{s},$
then
$$ k_{s}+ m_{s}-( k_{s}+1-t ) - \dim (( K \cap A_{s}) + (K^{\perp} \cap A_{s}))=-1, $$
implying
$t=0.$
\\
Now, consider the normal projective coordinate system for the space. Set
$ p_{j}= \displaystyle{\sum _{i=0}^{j} n_{i}}$
for 
$ 0 \leq j \leq r$
where
$n_{i}= m_{i}- m_{i+1}.$
Suppose that for each 
$i,$
$ 0 \leq i \leq k,$
the coordinate matrix of 
$ X_{i}$
is
$(\begin{array}{c c c c c c c c }
0 & \ldots & 0 & x_{p_{j-1}}^{i}& \ldots & x_{p_{j}-1}^{i} & \ldots & x_{n}^{i}
\end{array})^{T}$
where 
$ X_{i} \in A_{j} \setminus A_{j+1}$
for some
$ d \leq j \leq l.$
If
$ \mathsf{X}^{i}_{j}=(x_{p_{j-1}}^{i} \ldots x_{p_{j}-1}^{i})^{T},$
the coordinate matrix of 
$ H_{i}$
takes the form
$(\begin{array}{c c c c c c c }
b_{i}^{0} & \ldots & b_{i}^{p_{j-1}-1}& \alpha (\mathsf{ E}_{j}\mathsf{ X}^{i}_{j})^{T} & 0 & \ldots & 0
\end{array})$
for some real number
$ \alpha. $
Since
$H_{i}\nsupseteq A_{j} ,$
it is followed that
$ \alpha \neq 0.$
We can assume that
$ \alpha = 1.$
Assume that for each
$j,$
$ d \leq j \leq l,$
$ \mathsf{X}_{j}$
is an 
$ n_{j}$
by
$( k_{j}-k_{j+1})$
matrix
$  \begin{pmatrix}
\begin{array}{c|c|c}
\mathsf{X}^{k_{j+1}+1}_{j} & \ldots &  \mathsf{X}^{k_{j}}_{j}
\end{array}
\end{pmatrix}.$
Let
$\mathsf{K}$
and
$\mathsf{K}^{\perp}$
be
$$ 
\begin{pmatrix}
\begin{array}{c:c:c:c}
    & & & \bigzero_{p_{d-1} \times ( k_{d}-k_{d+1})} \\
   & \bigzero_{p_{l-2} \times ( k_{l-1}-k_{l})}  &  & \mathsf{X}_{d} \\
 \bigzero_{p_{l-1} \times ( k_{l}-k_{l+1})}  &  &  \reflectbox{$\ddots$} &  * \\
    & \mathsf{X}_{l-1} &  & \vdots \\
\mathsf{X}_{l} & * &  & * \\
* & *&  & * 
\end{array}\end{pmatrix}
$$
and
$$
\begin{pmatrix}
\begin{array}{c c c c c }
* & * & * \qquad \qquad (\mathsf{E}_{l} \mathsf{X}_{l})^{T}  & \qquad \bigzero_{( k_{l}-k_{l+1}) \times (n+1 - p_{l} )} &
\\
\hdashline
 \\
  * & * & (\mathsf{E}_{l-1} \mathsf{X}_{l-1})^{T} &  \bigzero_{( k_{l-1}-k_{l}) \times (n+1 - p_{l-1} )} &  
\\
  \hdashline
 &  &   \reflectbox{$\ddots$} \qquad \qquad \qquad & \qquad \qquad & \qquad \qquad
 \\
\hdashline
\\
* & \qquad  (\mathsf{E}_{d} \mathsf{X}_{d})^{T} & &  \bigzero_{( k_{d}-k_{d+1}) \times (n+1 - p_{d} )} \qquad \qquad \qquad \qquad & \qquad \qquad   
\end{array}
\end{pmatrix}
$$
respectively in which 
$(i+1)$-th
column of
$\mathsf{K}$
is the coordinate matrix of 
$ X_{i}$
and 
$(i+1)$-th
row of
$\mathsf{K^{\perp}}$
is the coordinate matrix of 
$ H_{i}.$
The matrix of the reflection in the pair
$ (K, K^{\perp})$
equals
$ \mathsf{I}_{n+1}-2 \mathsf{K}(\mathsf{K}^{\perp} \mathsf{K})^{-1} \mathsf{K}^{\perp}.$
The product
$\mathsf{K}^{\perp} \mathsf{K}$
and its inverse are the following block upper triangular matrices
$$
\mathsf{K}^{\perp} \mathsf{K}= \begin{pmatrix}
 \mathsf{X}_{l}^{T} \mathsf{E}_{l} \mathsf{X}_{l} & & * \\
   & \ddots  &   \\
\bigzero &  & \mathsf{X}_{d}^{T} \mathsf{E}_{d} \mathsf{X}_{d} \\
\end{pmatrix}
$$
and
$$
(\mathsf{K}^{\perp} \mathsf{K})^{-1}= \begin{pmatrix}
 (\mathsf{X}_{l}^{T} \mathsf{E}_{l} \mathsf{X}_{l})^{-1} & & * \\
   & \ddots  &   \\
\bigzero &  & (\mathsf{X}_{d}^{T} \mathsf{E}_{d} \mathsf{X}_{d})^{-1} \\
\end{pmatrix}.
$$
It can be seen that the matrices
$ \mathsf{K}(\mathsf{K}^{\perp} \mathsf{K})^{-1}$
and
$ \mathsf{K}(\mathsf{K}^{\perp} \mathsf{K})^{-1} \mathsf{K}^{\perp} $
are equal to
$$
 \begin{pmatrix}
 \begin{array}{c c c c }
 &  & \bigzero_{p_{d-1} \times (k+1)} \qquad \qquad \qquad \qquad & \\
 \hdashline
 \\
\qquad \quad \qquad & \qquad \bigzero_{n_{d} \times (k_{d+1}+1)} &  & \mathsf{X}_{d} (\mathsf{X}_{d}^{T} \mathsf{E}_{d} \mathsf{X}_{d})^{-1}\qquad \\
 \hdashline
 \\
 & \bigzero_{n_{d+1} \times (k_{d+2}+1)} \qquad \qquad  & \mathsf{X}_{d+1} (\mathsf{X}_{d+1}^{T} \mathsf{E}_{d+1} \mathsf{X}_{d+1})^{-1}  & * \\
 \hdashline
 & \qquad \qquad \qquad \qquad \reflectbox{$\ddots$} &  &  \\
 \hdashline
& \mathsf{X}_{l} (\mathsf{X}_{l}^{T} \mathsf{E}_{l} \mathsf{X}_{l})^{-1} \qquad \qquad * & \quad \ldots & * \\
\hdashline
\qquad * & * & \ldots & *
\end{array}
\end{pmatrix}
$$
and
$$
 \begin{pmatrix}
 \begin{array}{c c c c }
 &   & \bigzero_{p_{d-1} \times (n+1)} \qquad \qquad & \\
 \hdashline
 \\
 \ * & \mathsf{X}_{d} (\mathsf{X}_{d}^{T} \mathsf{E}_{d} \mathsf{X}_{d})^{-1} (\mathsf{E}_{d} \mathsf{X}_{d})^{T}  & \qquad \qquad \quad \bigzero_{n_{d} \times (n+1-p_{d})} &  \\
\hdashline
 & \qquad \qquad \qquad \qquad \ddots & & \\ 
\hdashline
\quad * &  \ldots \qquad \quad * & \mathsf{X}_{l} (\mathsf{X}_{l}^{T} \mathsf{E}_{l} \mathsf{X}_{l})^{-1} (\mathsf{E}_{l} \mathsf{X}_{l})^{-1} \qquad \qquad & \bigzero_{n_{l} \times (n+1-p_{l})} \\
 \hdashline
\quad * & \qquad \qquad \ldots &  \qquad * & \bigzero_{(n+1-p_{l}) \times (n+1-p_{l})} \qquad \\
\end{array}
\end{pmatrix}
$$
respectively. Therefore the matrix
$ \mathsf{I}_{n+1}-2 \mathsf{K}(\mathsf{K}^{\perp} \mathsf{K})^{-1} \mathsf{K}^{\perp}$
equals
$$
 \begin{pmatrix}
 \begin{array}{c c c c }
\ \mathsf{I}_{p_{d-1}}  &  & \bigzero_{p_{d-1} \times (n+1-p_{d-1})} \qquad \qquad &  \\
\hdashline
\\
*  & \mathsf{I}_{n_{d}}- 2 \mathsf{X}_{d} (\mathsf{X}_{d}^{T} \mathsf{E}_{d} \mathsf{X}_{d})^{-1} (\mathsf{E}_{d} \mathsf{X}_{d})^{T}   & \qquad \qquad \qquad \bigzero_{n_{d} \times (n+1-p_{d})} & \\
\hdashline
  & \qquad \qquad \qquad \qquad \qquad \ddots  &  &  \\
\hdashline
\quad * & \ldots  \qquad \qquad  * & \mathsf{I}_{n_{l}}-2 \mathsf{X}_{l} (\mathsf{X}_{l}^{T} \mathsf{E}_{l} \mathsf{X}_{l})^{-1} (\mathsf{E}_{l} \mathsf{X}_{l})^{T} \quad & \bigzero_{n_{l} \times (n+1-p_{l})} \qquad \\
\hdashline
\quad  * &  \qquad \qquad \ldots & * & \mathsf{I}_{n+1-p_{l}} \qquad \qquad \\
\end{array}
\end{pmatrix} 
$$
which is a block lower triangular matrix. It is easy to see that
$$ (\mathsf{I}_{n_{j}}-2 \mathsf{X}_{j} (\mathsf{X}_{j}^{T} \mathsf{E}_{j} \mathsf{X}_{j})^{-1} \mathsf{X}_{j} ^{T} \mathsf{E}_{j})^{T} \mathsf{E}_{j} (\mathsf{I}_{n_{j}}-2 \mathsf{X}_{j} (\mathsf{X}_{j}^{T} \mathsf{E}_{j} \mathsf{X}_{j})^{-1} \mathsf{X}_{j} ^{T} \mathsf{E}_{j})= \mathsf{E}_{j}$$
for 
$ d \leq j \leq l.$
\end{proof}
\label{t4}
\end{theorem}
\begin{theorem}
Every motion in a Cayley-Klein space of dimension 
$n$
is a composition of at most 
$n+1$
reflections in point-hyperplane pairs.
\begin{proof}
To prove the theorem we show that the matrix associated to a given motion is a product of at most 
$ n+1$
matrices of the form
$$  \mathsf{I}_{n+1}-2 \mathsf{X}( \mathsf{X}^{\perp} \mathsf{X})^{-1} \mathsf{X}^{\perp},$$
in which
$ \mathsf{X}$
and
$\mathsf{X}^{\perp}$
are the coordinate matrices of a point and a total polar of it respectively. In each Cayley-Klein space of dimension 
$n,$
with the exception of
$n=0,$
such a matrix represents a reflection in a point-hyperplane pair. For 
$n=0,$
the only motion of the space is the identity mapping which is a product of zero reflection. This motion is equal to
$ \mathsf{I}_{1}-2 \mathsf{X}( \mathsf{X}^{\perp} \mathsf{X})^{-1} \mathsf{X}^{\perp}$
for 
$\mathsf{X}=(1)$
and 
$ \mathsf{X}^{\perp}= (1).$
Suppose that the assertion is true in each Cayley-Klien space of dimension 
$n-1.$
Now, consider a Cayley-Klein space of dimension 
$n$
with the absolute figure
$$   \hat{Q_{0}} \supseteq \hat{A_{1}} \supseteq \hat{Q_{1}} \supseteq \hat{A_{2}} \supseteq \cdots \supseteq \hat{A_{r}} \supseteq \hat{Q_{r}} \supseteq \hat{A}_{r+1}= \varnothing.$$
Let 
$g$
be a motion of the space and 
$\mathsf{U}$
the matrix associated to it in the normal projective coordinate system. Thus 
$\mathsf{U}$
is a block lower triangular matrix such as
$$\begin{pmatrix}
 \mathsf{U}_{0} & & \bigzero \\
   & \ddots  &   \\
* &  & \mathsf{U}_{r} \\
\end{pmatrix}
$$
where
$ \mathsf{U}_{i} $
is a matrix of dimension 
$ n_{i}$
satisfying
$ \mathsf{U}_{i}^{T} \mathsf{E}_{i} \mathsf{U}_{i}= \mathsf{E}_{i}$
for
$0 \leq i \leq r.$
Let 
$X$
be the point with the coordinate matrix
$ \mathsf{X}=(\begin{array}{c c c c}
1 & 0 & \ldots & 0
\end{array})^{T}.$
If 
$g(X)=X,$
the first column of 
$ \mathsf{U} $
takes the form
$(\begin{array}{c c c c}
\alpha & 0 & \ldots & 0
\end{array})^{T}$
for some 
$ \alpha \neq 0.$
Since
$ \mathsf{U}_{0}^{T} \mathsf{E}_{0} \mathsf{U}_{0}= \mathsf{E}_{0},$
it is followed that
$$ \mathsf{U}=\begin{pmatrix}
\alpha & 0 & \ldots & 0& 0 \\
0 &  \mathsf{V}_{0} & & & \bigzero \\ 
0 &  & \mathsf{U}_{1} & &  \\
 \vdots & &  & \ddots  &   \\
0 & * & & & \mathsf{U}_{r}
\end{pmatrix}
$$
in which
$ \alpha^{2}=1$
and
 $ \mathsf{V}_{0}^{T} \mathsf{E}_{0}^{\prime} \mathsf{V}_{0}= \mathsf{E}_{0}^{\prime}$
for
$ \mathsf{E}_{0}= \begin{pmatrix} 
1 & 0 \\
0 & \mathsf{E}^{\prime}_{0}
\end{pmatrix}.$
It can be assumed that
$ \alpha =1.$
Consider a Cayley-Klein space of dimension
$n-1$
with the absolute figure
$$  \hat{Q_{0}^{\prime}} \supseteq \hat{A_{1}^{\prime}} \supseteq \hat{Q_{1}} \supseteq \hat{A_{2}} \supseteq \cdots \supseteq \hat{A_{r}} \supseteq \hat{Q_{r}} \supseteq \hat{A}_{r+1}= \varnothing.$$
where 
$ (\begin{array}{c c c }
x_{1} & \ldots & x_{n_{0}-1}
\end{array}) \mathsf{E}_{0}^{\prime} (\begin{array}{c c  c}
x_{1} & \ldots & x_{n_{0}-1}
\end{array})^{T} = 0$
is the equation of
$Q_{0}^{\prime}.$
By the induction hypothesis, for the matrix
$$ \mathsf{U}^{\prime}= \begin{pmatrix}
  \mathsf{V}_{0} & & & \bigzero \\ 
   & \mathsf{U}_{1} & &  \\
   & & \ddots  &   \\
 * & & & \mathsf{U}_{r}
\end{pmatrix}
$$
there exist matrices 
$ \mathsf{S}_{1}^{\prime}, \ldots, \mathsf{S}_{m}^{\prime} $
for 
$ m \leq n$
such that 
$\mathsf{U}^{\prime}= \mathsf{S}_{1}^{\prime} \cdots \mathsf{S}_{m}^{\prime}  $
in which
$$
\mathsf{S}_{i}^{\prime}= \mathsf{I}_{n}-2 \mathsf{X}_{i}^{\prime}( (\mathsf{X}_{i}^{\prime})^{\perp} \mathsf{X}_{i}^{\prime})^{-1} (\mathsf{X}_{i}^{\prime})^{\perp}
$$
for some point
$X_{i}^{\prime}$
and a total polar 
$(X_{i}^{\prime})^{\perp} $
of it with the coordinate matrices
$\mathsf{X}_{i}^{\prime}$
and
$ (\mathsf{X}_{i}^{\prime})^{\perp}$
respectively. It is easily seen that the hyperplane
$X_{i}^{\perp}$ 
with the coordinate matrix 
$ \mathsf{X}_{i}^{\perp}=( 0 \ (\mathsf{X}_{i}^{\prime})^{\perp}) $
is a total polar of the point
$X_{i}$
with the coordinate matrix
$ 
\begin{pmatrix}
0 \\
\mathsf{X}_{i}^{\prime}
\end{pmatrix}
$
in the 
$n$-dimensional
Cayley-Klein space. By setting
$$ \mathsf{S}_{i}= \mathsf{I}_{n+1}-2 \mathsf{X}_{i}( \mathsf{X}_{i}^{\perp} \mathsf{X}_{i})^{-1} \mathsf{X}_{i}^{\perp},$$
we get
$$
\mathsf{S}_{i}=\begin{pmatrix}
1 & 0 & \ldots & 0 \\
0 & & & \\
\vdots & & \mathsf{S}_{i}^{\prime} & \\
0 & & &
\end{pmatrix}
$$
for
$1 \leq i \leq m,$
implying
$ \mathsf{U}= \mathsf{S}_{1} \cdots \mathsf{S}_{m}.$
\par
If
$ g(X) \neq X, $
then
$ g(X)-X \notin Q_{0}$
or 
$ g(X)+X \notin Q_{0}$
in which
$g(X)\pm X$
is the point with the coordinate matrix
$\mathsf{U} \mathsf{X} \pm \mathsf{X}.$
In the former case, the unique total polar
$(g(X)-X)^{p}$
of
$g(X)-X$
contains
$g(X)+X$
and in the latter case we get
$ g(X)-X \in (g(X)+X)^{p}. $
Since
$g(X)$
is the harmonic conjugate of
$X$
with respect to
$g(X)-X$
and
$g(X)+X,$
it is followed that
$g(X)$
is mapped to
$X$
under the reflection 
$s$
in the pair 
$(g(X)-X, (g(X)-X)^{p})$
or
$(g(X)+X, (g(X)+X)^{p}).$
So 
$ (s \circ g)(X)=X$
and we get the result from the first case.
\end{proof}
\label{t5}
\end{theorem}
Therefore, in each Cayley-Klein space the set of all reflections in point-hyperplane pairs is a generator of the group of motions of the space. We show that the reflections in point-hyperplane pairs that both point and hyperplane are not regular are necessary for generating motions. Consider a Cayley-Klein space of dimension three with the absolute figure 
$$   \hat{Q_{0}} \supseteq \hat{A_{1}} \supseteq \hat{Q_{1}} \supseteq \hat{A_{2}} \supseteq \hat{Q_{2}} \supseteq \hat{A}_{3} \supseteq \hat{Q_{3}} \supseteq \varnothing, $$
in which each 
$Q_{i},$
$0 \leq i \leq 3,$
is a hyperplane of 
$A_{i}.$
All points of the space not in
$A_{1}$
and all planes not containing
$A_{3}$
are the only regular points and planes respectively. The total polar of each regular point is the plane 
$A_{1}$
and the total polar of each regular plane is the point
$A_{3}.$
In the normal projective coordinate system, the coordinate matrix of a regular point
$X$
and a regular plane
$H$
takes the form
$(\begin{array}{c c c c}
1 & x & y & z
\end{array})^{T}$
and 
$(\begin{array}{c c c c}
a & b & c & 1
\end{array})$
respectively. Also the coordinate matrices of the plane
$A_{1}$
and the point
$A_{3}$
are
$(\begin{array}{c c c c}
1 & 0 & 0 & 0
\end{array})$
and
$(\begin{array}{c c c c}
0 & 0 & 0 & 1
\end{array})^{T}$
respectively. So the matrices associated to the reflections in the pairs
$(X, A_{1})$
and
$(A_{3}, H)$
are
$$\begin{pmatrix}
-1 & 0 & 0 & 0 \\
-2x & 1 & 0 & 0 \\
-2y & 0 & 1 & 0 \\
-2z & 0 & 0 & 1 
\end{pmatrix}
\text{and} \ 
\begin{pmatrix}
1 & 0 & 0 & 0 \\
0 & 1 & 0 & 0 \\
0 & 0 & 1 & 0 \\
-2a & -2b & -2c & -1 
\end{pmatrix}.
$$
We show that the set of all such matrices can not span the group of motions of this Cayley-Klein space. For each matrix of the set the product of the second and the third diagonal entries is positive. Since each matrix is lower triangular, this positivity holds for every product of them. Let 
$X^{\prime}$
be a point of
$A_{1}\setminus A_{2}$
and 
$H^{\prime}$
a total polar of it different from
$A_{1}.$
The coordinate matrices of
$X^{\prime}$
and
$H^{\prime}$
are of the forms
$(\begin{array}{c c c c}
0 & 1 & y^{\prime} & z^{\prime}
\end{array})^{T}$
and
$(\begin{array}{c c c c}
a^{\prime} & 1 & 0 & 0
\end{array})$
respectively. The matrix associated to the reflection in the pair
$(X^{\prime},H^{\prime})$
is 
$$\begin{pmatrix}
1 & 0 & 0 & 0 \\
-2a^{\prime} & -1 & 0 & 0 \\
-2a^{\prime}y^{\prime} & -2y^{\prime} & 1 & 0 \\
-2a^{\prime}z^{\prime} & -2z^{\prime} & 0 & 1 
\end{pmatrix}.
$$
This matrix can not be a product of the above matrices.
 The authors report there are no competing interests to declare.

\end{document}